\documentclass[12pt]{amsart}

\usepackage[normalem]{ulem}

\usepackage{lineno}

\usepackage{soul}

\usepackage{graphicx}

\usepackage{amssymb, amsmath, amsthm, hyperref, euscript, color}

\usepackage[dvipsnames]{xcolor}

\usepackage{amscd}

\usepackage{tikz-cd}

\usepackage{bbm}

\usepackage{enumitem}

\usepackage[english]{babel}

\usepackage{mathtools}

\usepackage{makecell}

\numberwithin{equation}{section}

\theoremstyle{plain}

\newtheorem{theorem}{Theorem}[section]

\newtheorem{proposition}[theorem]{Proposition}

\newtheorem{lemma}[theorem]{Lemma}

\theoremstyle{definition}
\newtheorem{definition}[theorem]{Definition}
\newtheorem{remark}[theorem]{Remark}
\newtheorem{question}[theorem]{Question}

\newtheorem{example}[theorem]{Example}

\newtheorem{thm}{Theorem}

\makeatletter
\newsavebox\myboxA
\newsavebox\myboxB
\newlength\mylenA

\newcommand*\xoverline[2][0.75]{%
    \sbox{\myboxA}{$\m@th#2$}%
    \setbox\myboxB\null
    \ht\myboxB=\ht\myboxA%
    \dp\myboxB=\dp\myboxA%
    \wd\myboxB=#1\wd\myboxA
    \sbox\myboxB{$\m@th\overline{\copy\myboxB}$}
    \setlength\mylenA{\the\wd\myboxA}
    \addtolength\mylenA{-\the\wd\myboxB}%
    \ifdim\wd\myboxB<\wd\myboxA%
       \rlap{\hskip 0.5\mylenA\usebox\myboxB}{\usebox\myboxA}%
    \else
        \hskip -0.5\mylenA\rlap{\usebox\myboxA}{\hskip 0.5\mylenA\usebox\myboxB}%
    \fi}
\makeatother

\newcommand{\R}{\mathbb{R}}

\newcommand{\Z}{\mathbb{Z}}

\newcommand{\N}{\mathbb{N}}

\newcommand{\CP}{{\mathbb C\mkern-0.5mu\mathbb P}}
\newcommand{\CPbar}{\xoverline[0.75]{{\mathbb C\mkern-0.5mu\mathbb P}}}

\DeclareMathOperator{\id}{id}

\DeclareMathOperator{\Pd}{PD}

\DeclareMathOperator{\Kod}{Kod}

\newcommand{\simtimes}{\mathbin{\widetilde{\smash{\times}}}}

\DeclareMathOperator{\Sw}{SW}
\DeclareMathOperator{\cs}{\#}

\begin{document}

\author{Valentina Bais, Rafael Torres and Daniele Zuddas}

\title[Small exotic 4-manifolds with free abelian $\pi_1$]{Some examples of small irreducible exotic 4-manifolds with free abelian fundamental group}

\address{Scuola Internazionale Superiore di Studi Avanzati (SISSA), Via Bo\-nomea 265\\34136\\Trieste\\Italy}

\email{$\{vbais, rtorres\}$@sissa.it}

\address{Dipartamento di Matematica e Geoscienze\\Università degli Studi di Trieste \\ Via Valerio 12/b\\34127\\Trieste\\Italy.}

\email{dzuddas@units.it}

\subjclass[2020]{Primary 57R55; Secondary 57K40, 57R40}

\begin{abstract}We produce examples of pairwise non-diffeomorphic closed irreducible 4-manifolds with non-trivial free abelian fundamental group of rank less than three and small Euler characteristic. These exotic smooth structures become standard after taking a connected sum with a single copy of $S^2\times S^2$. The contributions of this paper include an explicit mechanism to computate the equivariant intersection form of 4-manifolds that are obtained via torus surgeries and a new stabilization result concerning exotic smooth structures with arbitrary fundamental group.

\end{abstract}

\maketitle

\section{Introduction.}

It has been noticed by several authors \cite{[AkhmedovPark2], [BaldridgeKirk], [Torres1], [Torres2]} that torus surgery based constructions of closed simply connected 4-manifolds of small Euler characteristic \cite{[AkhmedovPark1], [AkhmedovPark2], [AkhmedovBaykurPark], [BaldridgeKirk1], [BaldridgeKirk], [FintushelParkStern], [FintushelStern], [Szabo]} can be modified in order to produce small 4-manifolds with multiple choices of fundamental group, including free abelian  groups. The homeomorphism class of the 4-manifolds that have been obtained in this way, however, had not been previously determined. Merely knowing the isometry class of the intersection form over the integers of a non-simply connected 4-manifold does not suffice to determine its topological type. Instead, one must determine its equivariant intersection form up to isometry. Moreover, the 4-manifolds under consideration in this note are outside the stable range and Bass's cancellation result does not apply as in Hambleton-Teichner \cite[Corollary 3 (2)]{[HambletonTeichner]}; see Remark \ref{Remark Explanation}.

The main contribution of this note is to describe an algorithm to compute the equivariant intersection form of 4-manifolds that are constructed by using torus surgeries, and show that it is extended from the integers. This is of interest with respect to a conjecture due to Friedl-Hambleton-Melvin-Teichner \cite[Conjecture 1.3]{[FriedlhambletonMelvinTeichner]} (cf. \cite{[HambletonTeichner]}), which states that the equivariant intersection of any closed smooth 4-manifold with infinite cyclic fundamental group is extended from the integers.  Results of Freedman-Quinn \cite{[Freedman Quinn]}, Stong-Wang \cite{[StongWang]} and Hambleton-Kreck-Teichner \cite{[HKT]} are then used to pin down the homeomorphism type of several 4-manifolds with free abelian fundamental groups. 

Results of Gompf \cite{[Gompf1]} and Wall \cite{[Wall1]} say that two closed smooth oriented homeomorphic 4-manifolds become diffeomorphic after taking a connected sum with $n$ copies of $S^2\times S^2$ for some $n\in \N$. Baykur-Sunukjian \cite{[BaykurSunukjian]} showed that the value $n = 1$ suffices for any pair of closed smooth simply connected non-spin 4-manifolds that are related by a torus surgery. Another contribution of this note is to extend Baykur-Sunukjian's work on stabilization of simply connected 4-manifolds to arbitrary fundamental group under certain restrictions.

The following theorem is a sample of the production of these contributions.

\begin{thm}\label{Theorem Main} Let $k\in \{2, 3, 4, 5\}$. There are infinite collections\begin{center}$\{M_{p, k} : p\in \N\}$ and $\{N_{p, k} : p\in \N\}$\end{center}of closed smooth oriented pairwise non-diffeomorphic 4-manifolds that satisfy the following properties.
\begin{enumerate}
\item The 4-manifolds $M_{p, k}$ and $N_{p, k}$ are irreducible for all $p \neq 0$.
\item There are homeomorphisms 
\begin{align*}M_{p, k} &\cong_{\mathcal{C}^{0}} \cs_2\CP^2\cs_{k + 1}\CPbar^2\cs {(S^1\times S^3)} \\ 
N_{p, k} &\cong_{\mathcal{C}^{0}} \cs_{2}\CP^2\cs_{k + 1}\CPbar^2\cs {(T^2\times S^2)}\end{align*} for every $p\in \N$.

\item There are diffeomorphisms
\begin{align*} M_{p, k}\cs {(S^2\times S^2)}&\cong_{\mathcal{C}^{\infty}} M_{0, k}\cs {(S^2\times S^2)} \\ N_{p, k}\cs {(S^2\times S^2)}&\cong_{\mathcal{C}^{\infty}} N_{0, k}\cs {(S^2\times S^2)} \end{align*}for every $p \in \N$.
\end{enumerate}
Moreover, the 4-manifolds $M_{1, k}$ and $N_{1, k}$ admit a symplectic structure of symplectic Kodaira dimension two.
\end{thm}

\begin{remark}\label{Remark Explanation} A closed oriented 4-manifold $M$ with infinite cyclic fundamental group that satisfies the inequality $b_2(M) - |\sigma(M)|\geq 6$, where $b_2$ is the second Betti number and $\sigma$ is the signature, is within the stable range. It is known that the equivariant intersection form of such 4-manifolds is extended from the integers \cite[Corollary 3 (2)]{[HambletonTeichner]}. The topological prototypes of Theorem \ref{Theorem Main}, on the other hand, satisfy\begin{equation*}b_2(\cdot) - |\sigma(\cdot)| = 4.\end{equation*}
\end{remark}

\begin{remark}
    The 3-sphere $\varSigma^3\subset M_{p,k}$ that corresponds to the image of $\{\text{pt}\} \times S^3$ under the homeomorphism in item $(2)$ of Theorem \ref{Theorem Main} is locally flat embedded in $M_{p,k}$ and there is no continuous ambient isotopy that carries it onto a smoothly embedded 3-sphere for $p\neq 0$. Indeed, such a smoothly embedded 3-sphere $\varSigma^3$ would yield an $S^1\times S^3$ connected summand in $M_{p,k}$ that would contradict the irreducibility of our examples.
\end{remark}

The paper has been structured as follows. Section \ref{Section Preliminaries} contains a summary of the basic gadgets and a description of the main construction tool that we used. Section \ref{Section Equivariant Intersection} encompasses the definition of the equivariant intersection form, its properties, and a key example that is useful to compute it. The classification results of Freedman-Quinn, Stong-Wang and Hambleton-Kreck-Teichner are collected in Section \ref{Section Homeomorphism Results}. The definition of torus surgeries is given in Section \ref{Section Torus Surgeries}, where we also pin down the choices of framings, the notation that will be used as well as record proofs of properties of this surgical operation that are not found in the literature. Background results on symplectic geometry are given in Section \ref{Section Symplectic}. The mechanism and main result that are used to distinguish the diffeomorphism types of the examples constructed is described in Section \ref{SW}. A characterization of the second homotopy group of 4-manifolds with certain free abelian fundamental groups is given in Section \ref{str}. The key component of the procedure to compute the equivariant intersection forms of 4-manifolds that are obtained via torus surgeries is detailed in Section \ref{construction}; another method to compute intersections of surfaces is outlined in \cite{[CollinsPowell]}. Section \ref{Section BK Construction} contains a description of a construction of small symplectic 4-manifolds due to Baldridge-Kirk \cite{[BaldridgeKirk]}, where we also investigate further the produce of their fundamental group computations. The examples of Theorem \ref{Theorem Main} arise from Baldridge-Kirk's work and their construction is described in Section \ref{Z}. Their homeomorphism classes are identified in Section \ref{Section Infinite Families} by computing the equivariant intersection form through the intersection of the lifts of surfaces to the corresponding universal cover. A procedure to obtain such surfaces in 4-manifolds that are obtained from torus surgeries is described in Section \ref{construction}. Our result on stabilizations with a single copy of $S^2\times S^2$ is presented in Section \ref{stab}. A proof of Theorem \ref{Theorem Main} is found in Section \ref{Section Proof}.

Throughout the paper, we will always work with smooth oriented manifolds. Homeomorphism and diffeomorphism are assumed to be orientation preserving and homology groups are taken with integer coefficients.

\subsection{Acknowledgements}
We would like to thank Daniel Kasprowski for useful email exchanges on the structure of the second homotopy group of closed orientable 4-manifolds. We thank Mark Powell for his valuable input and comments on an earlier version, which helped us improve the paper. The first author is sincerely grateful to Younes Benyahia and Oliviero Malech for the useful discussions and to Emanuele Pavia for the help in computing spectral sequences and for the enjoyable work at Barcola. The first and the third authors are members of GNSAGA – Istituto Nazionale di Alta Matematica ‘Francesco Severi’, Italy. The first author was partially supported by GNSAGA. 

\section{Preliminary results and construction tools.}\label{Section Preliminaries}

\subsection{Equivariant intersection form}\label{Section Equivariant Intersection} In the following we recall the definition of the equivariant intersection form of a closed oriented 4-manifold $X$, since this will be a crucial ingredient for our purposes in light of Theorem \ref{classification}. The group ring $\mathbb{Z}[\pi_1(X)]$ is endowed with the involution obtained by extending linearly the inverse map automorphism of $\pi_1(X)$ and it acts on $H_2(\widetilde{X})$ via Deck transformations, making it a $\Z[\pi_1(X)]$-module. Here $\widetilde{X}$ denotes the universal cover of $X$. We will always suppose that a base-point is fixed, although it will be omitted from our notation if it is not explicitly needed.

\begin{definition}\label{defmain}
    The equivariant intersection form of a closed oriented 4-manifold $X$ is the hermitian form
    $$\widetilde{I}_X: H_2(\widetilde{X}) \times H_2(\widetilde{X}) \to \mathbb{Z}[\pi_1(X)]$$
    defined by setting
    \begin{equation}\label{rule}
        \widetilde{I}_X(a,b)\coloneq \sum_{g \in \pi_1(X)} I_{\widetilde{X}}(a,g_*(b)) g^{-1}
    \end{equation}
    for any $a,b \in H_2(\widetilde{X})$, where $I_{\widetilde{X}}$ denotes the intersection form of $\widetilde{X}$ over the integers and $g_*$ is the image of $g \in \pi_1(X)$ under the identification with the group of Deck transformations of the universal cover. 

    The radical of the equivariant intersection form is the subset
    \[R(\widetilde{I}_X)\coloneq \{a \in H_2(\widetilde{X}) \mid \widetilde{I}_X(a,b)=0 \ \forall b \in H_2(\widetilde{X}) \}.\]
    
    In the following, the identification $H_2(\widetilde{X}) \cong \pi_2(X)$ induced by the Hurewicz isomorphism is understood.
\end{definition}

For our purposes, we are interested in the study of the equivariant intersection form of a 4-manifold up to the following equivalence relation.

\begin{definition}\label{isometry}
    Let $X$ and $Y$ be two closed oriented 4-manifolds. An isometry between their equivariant intersection forms $\widetilde I_X$ and $\widetilde I_{Y}$ is a pair $(g,h)$, where $$g: \pi_1(X) \to \pi_1(Y)$$ is an isomorphism of fundamental groups and $$h: H_2(\widetilde{X}) \to H_2(\widetilde{Y})$$ is $g$-invariant and satisfies $$\widetilde{I}_X(a,b)=\widetilde{I}_{Y}(h(a),h(b))$$ for all $a,b \in H_2(\widetilde{X}).$

    If such an isometry exists, we will say that the two equivariant intersection forms are isometric and write $\widetilde I_X\cong \widetilde I_{Y}$.
\end{definition}

\begin{remark}\label{reduced}
        The hermitian forms induced by the equivariant intersection forms of two closed oriented 4-manifolds $X$ and $Y$ on the quotients $H_2(\widetilde{X})/R(\widetilde{I}_X)$ and $H_2(\widetilde{Y})/R(\widetilde{I}_{Y})$ are non-singular (provided $H^3(\pi_1(X), \pi_2(X)) = 0 = H^3(\pi_1(Y), \pi_2(Y))$) and they are isometric if and only if $\widetilde{I}_X \cong \widetilde{I}_{Y}$ (see \cite[Remark 1.4]{[HKT]}). In particular, this implies that there is no loss in generality in studying the isometry class of the equivariant intersection form of a 4-manifold over the quotient $H_2(\widetilde{X})/R(\widetilde{I}_X)$.
    \end{remark}

The topological prototypes of Theorem \ref{Theorem Main} can be realized as a connected sum $X = M\cs X_0$, where $M$ is a closed simply connected 4-manifold and $X_0$ is either $S^1\times S^3$ or $T^2\times S^2$. The equivariant intersection form of such a connected sum is said to be extended from the integers. In the following, we recall the general definition \cite{[Hambleton]}. 

    \begin{definition}\label{extendedint}
    The equivariant intersection form $\widetilde{I}_X$ of $X$ is extended from the integers if $$\widetilde{I}_X=Q\otimes_{\mathbb{Z}} \mathbb{Z}[\pi_1(X)]$$  for some integral unimodular form $Q$.
\end{definition}

\begin{remark}\label{extended}
         Let $U \subset \widetilde{X}$ be an open subset with the property that the universal covering map $\widetilde {p}: \widetilde{X} \to X$ sends it homeomorphically onto an open subset $\widetilde{p}(U)\subset X$. If $U$ contains a collection of embedded surfaces generating the quotient $H_2(\widetilde{X})/R(\widetilde{I}_X)$ as a $\Z[\pi_1(X)]$-module, then $\widetilde{I}_X$ is extended from the integers \ref{extendedint}. 
    \end{remark}

%

\begin{example}\label{example}
Let $X$ be a closed smooth oriented 4-manifold with fundamental group $\pi_1(X) \cong \Z$ and a fixed handlebody decomposition. Without loss of generality, we can assume that just one among the 1-handles gives a non-trivial generator of $\pi_1(X)$. Let $X^2$ be the union of all handles up to index 2 and $h^1$ the non-trivial 1-handle in $\pi_1(X)$. From the recipe for building handlebody decompositions of cyclic coverings (see e.g. \cite[Section 6.3]{[GompfStipsicz]}) it follows that the pre-image of the open subset $V:=\text{Int}(X^2 \setminus h^1)\subset X$ via the universal covering map is a disjoint union of copies of $V$. By Remark \ref{extended}, the equivariant intersection form of $X$ is extended from the integers and coincides with the intersection form $I_X$ of $X$ whenever the pre-image of $V\subset X$ under the universal covering map is $H_2(\widetilde{X})$-surjective. A similar argument can be applied to the case in which $\pi_1(X)\cong \Z^2$. 
\end{example}

\subsection{Classification up to homeomorphism}\label{Section Homeomorphism Results} Results of Freedman-Quinn \cite{[Freedman Quinn]}, Stong-Wang \cite{[StongWang]} and Hambleton-Kreck-Teichner \cite{[HKT]} indicate that the equivariant intersection form encodes enough topological information of a closed smooth 4-manifold in order to pin down its homeomorphism class under some restrictions of the fundamental group. We summarize their results in the following statement. Recall that an oriented 4-manifold $X$ has\begin{enumerate}\item type (I) if its universal cover is non-spin, \item type (II) if it is spin and \item type (III) if it is non-spin but with a spin universal cover.\end{enumerate}

\begin{theorem}\label{classification}Let $X$ and $Y$ be closed smooth oriented 4-manifolds.

$\bullet$ Freedmann-Quinn \cite{[Freedman Quinn]}, Stong-Wang \cite{[StongWang]}. If the fundamental group of $X$ is infinite cyclic and there is an isometry $h: \widetilde{I}_X\rightarrow \widetilde{I}_{Y}$ between the equivariant intersection forms of $X$ and $Y$ as in Definition \ref{isometry}, then $h$ is realized by a homeomorphism $X\rightarrow Y$.

$\bullet$ Hambleton-Kreck-Teichner \cite[Theorem A]{[HKT]}. If the fundamental group of $X$ is free abelian of rank two and there is an isometry $h: \widetilde{I}_X\rightarrow \widetilde{I}_{Y}$ between the equivariant intersection forms of $X$ and $Y$ as in Definition \ref{isometry}, then $h$ is realized by a homeomorphism $X\rightarrow Y$ provided $X$ and $Y$ are of the same type. 
\end{theorem}

A proof of the first clause of Theorem \ref{classification} under the assumption of the existence of a homotopy equivalence $X\rightarrow Y$ can be found in \cite[Theorem 2.1]{[Powell]}.


\subsection{Torus surgeries and framings}\label{Section Torus Surgeries}Let $T\subset X$ be a smoothly embedded framed 2-torus with trivial normal bundle. Let $a$ and $b$ be simple loops in $T$ whose homotopy classes generate $\pi_1(T) = \Z\times \Z$ so that any primitive curve $\gamma_T\subset T$ can be expressed as $\gamma_T = a^r b^s$ for some relatively prime $r, s\in \Z$. Let\begin{equation}\label{Framing}\phi: \nu(T)\rightarrow T^2\times D^2\end{equation} be a framing of $T$ and let $\pi: T^2\times D^2\rightarrow T^2$ be the projection onto the first factor. Here $\nu(\cdot)$ denotes a tubular neighbourhood of the submanifold in brackets. The framing (\ref{Framing}) provides an isomorphism $H_1(\partial \nu(T)) \cong H_1(T^2)\oplus \Z$, where the $\Z$-factor is generated by the isotopy class of the meridian $\mu_T$ of the 2-torus $T$. For $d\in \partial D^2$, we set
\begin{align*}
S^1_a &= \phi^{-1}(\pi(a)\times \{d\})\\
S^1_b &= \phi^{-1}(\pi(b)\times \{d\})
\end{align*}
for the push-offs of the curves $a$ and $b$, respectively, in $\partial(X\setminus\nu(T))$ under (\ref{Framing}). In $\nu(T)$, the loop $S^1_a$ is homologous to $a$ and the loop $S^1_b$ is homologous to $b$. Notice that the triple of loops $\{S^1_a, S^1_b, \mu_T\}$ generates the first homology group $H_1(\partial \nu(T))$ \cite[Section 3]{[FintushelStern]}. The surgery curve $\gamma\subset \partial(X\setminus \nu(T))$ is the image of a primitive curve $\gamma_T\subset T$ under the framing (\ref{Framing}). The choice of surgery curve $\gamma = S^1_b$ suffices for the purposes of this note. The choice of framing used in the sequel is one of the following two canonical choices.

$\bullet$ If the 4-manifold $X$ has a symplectic structure for which the 2-torus $T\subset X$ is Lagrangian, we will choose (\ref{Framing}) to be the Lagrangian framing \cite{[HoLi]}.

$\bullet$ If\begin{equation}\label{Homology Class}[T] = 0 \in H_2(X),\end{equation} there is preferred framing known as the null-homologous framing \cite[Section 5.3]{[HoLi]}. Under the null-homologous framing of $T$, we have that\begin{equation}\label{Nullhomologous framing homology}[\gamma] = 0 \in H_1(X\setminus \nu(T))\end{equation} for the surgery curve. 

For any pair of co-prime integers $p,q \in \Z$, we construct a smooth 4-manifold\begin{equation}\label{Torus Surgery}X_{T, \gamma}(p/q) = (X\setminus \nu(T)) \cup_{\varphi} (T^2\times D^2),\end{equation}where the diffeomorphism $\varphi: T^2\times \partial D^2 \rightarrow \partial (X\setminus \nu(T))$ satisfies\begin{equation}\label{Gluing Map Induced}\varphi_{\ast}([\{pt\}\times \partial D^2]) = q[\gamma] + p[\mu_T]\in H_1(\partial (X\setminus \nu(T))).\end{equation} This cut-and-paste operation is known as a $(p/q)$-torus surgery on $T$ along $\gamma$. The submanifold \begin{equation}\label{Core Torus}\widehat{T}_{p/q}:= T^2\times \{(0, 0)\}\subset T^2\times D^2\end{equation}of $X_{T, \gamma}(p/q)$ is known as the core 2-torus of the surgery. Notice that\begin{equation}X\setminus \nu(T) = X_{T, \gamma}(p/q)\setminus \nu(\widehat{T}_{p/q})\end{equation}due to the local nature of the construction (\ref{Torus Surgery}). In the sequel, we will use the notation $X_{T, \gamma}(p) = X_{T, \gamma}(p/1)$ and call such a surgery a $p$-torus surgery. For further details on this cut-and-paste operation, we refer the reader to  \cite[Section 2.1]{[BaldridgeKirk]}, \cite{[HoLi]}, \cite[p. 310]{[GompfStipsicz]}.

The following proposition summarizes several basic topological properties of this cut-and-paste operation.

\begin{proposition}\label{surgery}
    The signature and the Euler characteristic of a closed oriented 4-manifold $X$ are invariant under torus surgeries. The fundamental group of a 4-manifold obtained by $(p/q)$-torus surgery is given by\begin{equation*}\pi_1(X_{T, \gamma}(p/q)) = \frac{\pi_1(X\setminus \nu(T))}{\langle [\mu_T]^p[\gamma]^q = 1\rangle},\end{equation*}where $\langle [\mu_T]^p[\gamma]^q = 1\rangle$ is the normal subgroup generated by $[\mu_T]^p[\gamma]^q$.
    
    Let $T \subset X$ be a homologically essential embedded 2-torus with trivial normal bundle and $\gamma_T \subset T$ a curve whose homology class is primitive in $H_1(X)$. If $q \neq 0$, the following facts hold.
    \begin{itemize}
        \item  The homology class of the surgery curve $\gamma \subset \partial (X \setminus \nu(T))$ given by the push-off of $\gamma_T \subset T$ in $X_{T, \gamma}(p/q)$ has order $q$ in $H_1(X_{T, \gamma}(p/q))$.
        \item The first Betti numbers satisfy the relation\begin{center}$b_1(X_{T, \gamma}(p/q))=b_1(X)-1$.\end{center} 
        \item The second Betti numbers satisfy the relation\begin{center}$b_2(X_{T, \gamma}(p/q))=b_2(X)-2$.\end{center}
        Furthermore, if $q=1$, the core 2-torus $\widehat{T}_{p} \subset X_{T, \gamma}(p)$ is null-homologous. 
    \end{itemize}
\end{proposition}

Torus surgeries can be reversed in the following sense: the original 4-manifold $X$ is obtained from $X_{T, \gamma}(p/q)$ by surgery on the core 2-torus $\widehat{T}_{p/q}$. A detailed description of this procedure for $p$-torus surgeries is given in the following lemma, which will be used in Section \ref{stab}. For the convenience of the reader, we provide two proofs of this fact. 

\begin{lemma}\label{inverting} Let $X_{T, \gamma}(p)$ be the 4-manifold defined in (\ref{Torus Surgery}). The core 2-torus $\widehat T_{p}$ is framed by identifying one of its tubular neighbourhoods with the copy of $T^2 \times D^2$ which is glued to $X \setminus \nu(T)$ during the surgery procedure. The 4-manifold $X$ is the result of $0$-torus surgery on the core 2-torus $\widehat T_{p} \subset X_{T, \gamma}(p)$ along a curve whose push-off corresponds to the meridian $\mu_T$ of the surgery 2-torus $T\subset X$ under the identification
\begin{center} $X_{T, \gamma}(p) \setminus \nu(\widehat{T}_{p})=X \setminus \nu(T).$\end{center} 

\end{lemma}

\begin{proof}[First proof of Lemma \ref{inverting}]
     Fix a framing $\phi: \nu(T)\rightarrow T^2\times D^2$ of the surgery 2-torus $T\subset X$. Without loss of generality, we can assume that the surgery curve $\gamma$ corresponds to the second $S^1$-factor. Following \cite[Section 6.3]{[Akbulutbook]}, we choose the gluing map $\varphi: T^2\times \partial D^2 \rightarrow \partial (X\setminus \nu(T))$ defining $X_{T, \gamma}(p)$ (\ref{Torus Surgery}) in such a way that the induced map in first homology is given by the matrix
    \begin{equation}\varphi_*=
         \begin{pmatrix}\label{martix}
1 & 0 & 0 \\
0 & 0 & 1\\
0 & -1 & p\\
\end{pmatrix} 
    \end{equation}
with respect to the basis $\{[S^1\times \{*\} \times \{*\}], [\{*\} \times S^1 \times \{*\}], [\{*\} \times \{*\} \times \partial D^2]\}$ of $H_1(S^1 \times S^1 \times \partial D^2)$, where the 2-torus $T^2$ is written as the product of two $S^1$-factors. 

 The identification
\begin{equation}
    X=(X \setminus \nu(T)) \cup_{\text{Id}} \nu(T) \cong (X_{T, \gamma}(p) \setminus \nu(\widehat T_{p})) \cup_{\psi} (T^2 \times D^2)
\end{equation}
shows that $X$ is obtained by surgery on $X_{T, \gamma}(p)$ along the core 2-torus $\widehat T_{p}$, where the gluing map
\begin{equation}
    \psi:T^2 \times \partial D^2 \to \partial (X_{T, \gamma}(p) \setminus \nu(\widehat T_{p}))
\end{equation}
sends the curve $\{*\} \times \partial D^2$ to a loop $\alpha$ which  under the identification  $\partial (X_{T, \gamma}(p) \setminus \nu(\widehat T_{p}))= \partial( X \setminus \nu( T))$ coincides with the meridian $\mu_T$ of the surgery 2-torus $T \subset X$. Notice that $\alpha$ is also the push-off of an essential curve $\alpha_{\widehat T_{p}}$ on the core 2-torus $\widehat T_{p}$. Indeed, by inverting the matrix (\ref{martix}), one can check that the homology class of $\mu_T\subset \partial(X \setminus \nu(T))$ is mapped to the one defined by the second $S^1$-factor of $S^1 \times S^1 \times D^2$, which in turn is identified with $\nu(\widehat{T}_{p}) \subset X_{T, \gamma}(p)$ by our choice of the framing. In particular, the loop $\alpha_{\widehat T_{p}}\subset \widehat T_{p}$ giving $\alpha$ as a push-off corresponds to the second $S^1$-factor under the identification $\widehat T_{p} \cong S^1 \times S^1 \times \{0\} \subset S^1 \times S^1 \times D^2$. Since $\psi_*([\{*\} \times \partial D^2])=[\alpha]$ the conclusion follows.
    \end{proof}

Our second proof of Lemma \ref{inverting} employs the useful technology of round 2-handles.

\begin{proof}[Second Proof of Lemma \ref{inverting}]
     A result of Baykur-Sunukjian \cite[Lemma 1]{[BaykurSunukjian]} implies that there is a cobordism\begin{equation} V = (I\times X) \cup R_2\end{equation}
     such that $\partial_-V = X$ and $\partial_+ V = X_{T, \gamma}(p)$, where $R_2 = S^1\times D^2\times D^2$ is a five-dimensional round 2-handle glued to $\partial_+(I \times X)$ via a map
     \begin{equation}\label{glue}
         \phi: R_2= S^1 \times D^2 \times D^2 \to \partial_+(I \times X)
     \end{equation}
     sending $S^1 \times \partial D^2 \times \{0\}$, $S^1 \times \partial D^2 \times D^2$ and $\{*\}\times \partial D^2 \times \{0\}$ homeomorphically onto $\{1\} \times T$, $\{1\} \times \nu(T)$ and $\{1 \} \times \gamma$ respectively. In particular, the curve $\{*\} \times \{*\} \times \partial D^2\subset S^1 \times \partial D^2 \times D^2$ is identified with a meridian of the surgery 2-torus $T \subset X$. Turning the cobordism $V$ upside down, we get a cobordism
     \begin{equation}\label{Cobordism 3}
     W = (I\times X_{T, \gamma}(p)) \cup R'_2\end{equation}
     where $R_2'=R_2$ is glued to $\partial_+(I \times X_{T, \gamma}(p))$ via a map
     \begin{equation}\label{glue2}
         \psi: S^1 \times D^2 \times \partial D^2 \to \partial_+(I \times X_{T, \gamma}(p))
     \end{equation}
      sending $S^1 \times \{0\} \times \partial D^2$, $S^1 \times D^2 \times \partial D^2$ and $\{*\}\times \{0\} \times \partial D^2$ homeomorphically onto $\{1\} \times \widehat{T}_{p}$, $\{1\} \times \nu(\widehat{T}_{p})$ and $\{1 \} \times \alpha$ respectively, where $\alpha \subset \widehat{T}_{p}$ is an essential curve. \cite[Lemma 1]{[BaykurSunukjian]} by Baykur-Sunukjian implies that $X$ is obtained from $X_{T,\gamma}(p)$ via a torus surgery on the core 2-torus $\widehat{T}_{p}$ along $\alpha$. Since the image of $\{*\} \times \{*\} \times \partial D^2 \subset S^1 \times \partial D^2 \times D^2$ under (\ref{glue2}) is a push-off of $\alpha$ which has been identified with a meridian of $T$ by (\ref{glue}), the result follows.
\end{proof}

\subsection{Luttinger surgery, irreducibility, and the symplectic Kodaira dimension}\label{Section Symplectic} We survey in this section several results of Auroux-Donaldson-Katzarkov \cite{[AurouxDonaldsonKatzarkov]}, Kotschick \cite{[Kotschick]}, Ho-Li \cite{[HoLi]} and Luttinger \cite{[Luttinger]} that are used to conclude on the properties of the symplectic 4-manifolds of Theorem \ref{Theorem Main}, following the notation of Section \ref{Section Torus Surgeries}.  Recall that a 4-manifold $X$ is said to be minimal if it does not contain any smoothly embedded 2-spheres of self-intersection $-1$. A $(\pm 1/q)$-torus surgery on a symplectic 4-manifold along a Lagrangian 2-torus with respect to the Lagrangian framing is called a Luttinger surgery \cite{[AurouxDonaldsonKatzarkov], [Luttinger]}. 

\begin{theorem}\label{Theorem Luttinger Surgery}Let $(X, \omega)$  be a symplectic 4-manifold that contains a Lagrangian 2-torus $T\subset X$ of self-intersection 0 equipped with the Lagrangian framing.\

$\bullet$ (Luttinger \cite{[Luttinger]}, Auroux-Donaldson-Katzarkov \cite{[AurouxDonaldsonKatzarkov]}).  The 4-manifold $X_{T, \gamma}(\pm 1/q)$ admits a symplectic structure.

$\bullet$ (Ho-Li \cite[Proposition 3.1]{[HoLi]}). If $(X, \omega)$ is minimal, then $X_{T, \gamma}(\pm 1/q)$ is minimal.

\end{theorem}

Let $(X, \omega)$ be a minimal symplectic four-manifold whose canonical class is $K_\omega=c_1(X,J)$, where $J$ is an $\omega$-compatible almost-complex structure on $X$.

\begin{definition}[The symplectic Kodaira dimension]\label{Definition Kodaira} Li \cite[Definition 2.2]{[HoLi]}. The symplectic Kodaira dimension $\Kod(X, \omega)$ of a minimal symplectic four-manifold $(X, \omega)$ is defined by\[
\Kod(X, \omega) =
\begin{cases}
-\infty & \text{if}\enspace K_\omega\cdot [\omega] < 0 \; \text{or} \enspace K_\omega \cdot K_\omega < 0,\\
0 & \text{if}\enspace K_\omega\cdot [\omega] = 0 \enspace \text{and} \enspace K_\omega \cdot K_\omega = 0,\\
1 & \text{if} \enspace K_\omega\cdot [\omega] > 0\enspace  \text{and} \enspace K_\omega \cdot K_\omega = 0 \enspace,\\
2 & \text{if}\enspace K_\omega\cdot [\omega] > 0 \enspace \text{and} \enspace K_\omega \cdot K_\omega > 0. 
\end{cases}
\]

The symplectic Kodaira dimension of a non-minimal symplectic 4-manifold is defined as the symplectic Kodaira dimension of any of its minimal models. 
\end{definition}

Ho-Li's study of the effect that Luttinger surgery has on the symplectic canonical class $K_{\omega}$ and on the symplectic class $[\omega]$ of a symplectic 4-manifold $(X, \omega)$ yielded the following result. 

\begin{theorem}[Ho-Li \cite{[HoLi]}]\label{Theorem Luttinger Surgery Kodaira} The symplectic Kodaira dimension is preserved under Luttinger surgery.

\end{theorem}

The reader is directed towards \cite[Section 3]{[HoLi]} for further details. Besides minimality, we are interested in 4-manifolds with the following property.

\begin{definition}
   A smooth 4-manifold $X$ is irreducible if for every smooth connected sum decomposition $X = X_1\cs X_2$ either $X_1$ or $X_2$ is a homotopy 4-sphere \cite[Definition 10.1.17]{[GompfStipsicz]}. 
\end{definition}

In order to conclude on the irreducibility of our examples, we will employ the following result.

\begin{theorem}[Kotschick \cite{[Kotschick]}]\label{Theorem Irreducibility} Every minimal symplectic 4-manifold $(X, \omega)$ with residually finite fundamental group and $b_2^+(X) > 1$ is irreducible. 
\end{theorem}

\begin{remark}
    Free abelian groups are residually finite.
\end{remark}

\subsection{Infinitely many diffeomorphism types and irreducibility}\label{SW}We make use of work of Fintushel-Park-Stern \cite{[FintushelParkStern]} based on the Morgan-Mrowka-Szab\'o formula \cite{[MorganMrowkaSzabo]} to distinguish the diffeomorphism types of the 4-manifolds constructed in this paper via their Seiberg-Witten invariants \cite[Section 2]{[FintushelParkStern]}. We regard the Seiberg-Witten invariant of a closed smooth oriented 4-manifold $X$ as a function\begin{center}$\Sw_X:\{k\in H_2(X): \Pd(k)\equiv w_2(X) \mod 2\}\rightarrow \Z$,\end{center}where $\Pd(k)$ denotes the Poincar\'e dual of $k$; see \cite{[FintushelParkStern], [AkhmedovBaykurPark]} for further details. The following result summarizes the results of Fintushel-Park-Stern \cite[Theorem 1, Corollary 1]{[FintushelParkStern]} that we will need in the sequel. We use the notation introduced in Section \ref{Section Torus Surgeries}. 

\begin{theorem}[Fintushel-Park-Stern \cite{[FintushelParkStern]}]\label{Theorem Smooth Structures} Let $X'$ be a closed smooth 4-manifold that contains a smoothly embedded null-homologous 2-torus $T'\subset X'$ and let $\gamma'\subset T'$ be a simple loop such that its push-off is null-homologous in $X'\setminus \nu(T')$. Furthermore, suppose that the 4-manifold $X'_{T', \gamma'}(0)$ has exactly one Seiberg-Witten basic class up to sign. The collection $\{X'_{T', \gamma'}(1/p): p\in \N\}$ consists of pairwise non-diffeomorphic elements. 
\end{theorem}


The gauge-theoretical background on the usage of Theorem \ref{Theorem Smooth Structures} has already been well-documented in the literature \cite{[AkhmedovBaykurPark], [BaldridgeKirk], [FintushelParkStern], [FintushelStern],[Torres2]}. The mechanism to apply Theorem \ref{Theorem Smooth Structures} to our examples is summarized in the diagram\[\begin{tikzcd}[row sep=4em]
X= X'_{T', \gamma'}(0) \arrow{r}{(1)} \arrow[swap]{d}{(2)} & X_{T, \gamma}(1) = X' \arrow{d}{(3)} \\%
X_{T, \gamma}(p) \arrow{r}{\id}  & X'_{T', \gamma'}(1/(p - 1))
\end{tikzcd}\]where the numbered arrows indicate the following steps in the construction; cf. Proposition \ref{surgery}.

(1).\ The irreducible symplectic 4-manifold $X$ on the left of the diagram contains a homologically essential Lagrangian 2-torus $T$ of self-intersection zero. Using the Lagrangian framing, perform a $(\pm 1)$-torus surgery on $T\subset X$ along a given surgery curve $\gamma\subset T$ that is primitive in $H_1(X)$ to obtain the symplectic 4-manifold $X_{T, \gamma}(\pm 1)$ on the upper right side of the diagram. Set $X':= X_{T, \gamma}(\pm 1)$ and let $T'\subset X'$ be the null-homologous core 2-torus of the surgery; see Proposition \ref{surgery}. Notice that $X= X'_{T', \gamma'}(0)$ as it was proven in Lemma \ref{inverting}.  A result of Taubes \cite{[Taubes]} says that the canonical classes $k = - c_1(X, J)$ and $k' = - c_1(X', J')$ are basic classes of $X$ and $X'$, respectively, and their Seiberg-Witten invariants are $\Sw_{X}(\pm k) = \pm 1 = \Sw_{X'}(\pm k')$. An argument due to Akhmedov-Baykur-Park implies that, under the circumstances considered in our construction, these are the only basic classes of $X$ and $X'$ \cite[Section 4.1]{[AkhmedovBaykurPark]}.

(2).\ Perform a $p$-torus surgery on $T\subset X$ along $\gamma\subset X$ using the Lagrangian framing to obtain a 4-manifold $X_{T, \gamma}(p)$.

(3).\ Perform a $(1/ (p - 1))$-torus surgery on $T'\subset X'$ along the surgery curve $\gamma'\subset T'$ using the null-homologous framing to produce the 4-manifold $X'_{T'. \gamma'}(p)$. Notice that \begin{equation*}X_{T, \gamma}(p) = X'_{T', \gamma'}(1/(p - 1)).\end{equation*} The Morgan-Mrowka-Szab\'o formula yields
\begin{equation*}\label{MMS Formula} \begin{split}\Sw_{X_{T', \gamma'}(1/(p - 1))}(k_{1/p}) &= \Sw_{X'}(k') + (p - 1)\cdot \Sw_{X}{(k)}\\
&= 1 + (p - 1) = p \end{split}\end{equation*}
where the basic class $k_{1/p}$ of $X_{1/(p - 1)}$ corresponds to the basic class $k_X$ of $X$. Thus, the collection $\{X_{T, \gamma}(p): p\in \N\}$ consists of pairwise non-diffeomorphic  4-manifolds given that their Seiberg-Witten invariants are pairwise different. 


\begin{remark}\label{Remark Szabo Computation}We point out that an argument due to Szab\'o \cite{[Szabo]} that uses a special case of the Morgan-Mrowka-Szab\'o formula \cite[Theorem 3.4]{[Szabo]} can also be used to reach the same conclusion on the Seiberg-Witten invariant of the collection $\{X_{T, \gamma}(p): p\in \N\}$.

\end{remark}


The following result will allow us to conclude that the collection $\{X_{T, \gamma}(p): p\in \N\}$ is made of irreducible 4-manifolds. It generalizes a result of Park on simply connected 4-manifolds \cite[Theorem 1]{[Park]} and its proof follows Park's argument to a great degree.

\begin{theorem}\label{Theorem Irreducible Diffeo}Let $Z$ be a closed minimal symplectic 4-manifold with $b_2^+(Z) > 1$ and let $M$ be a closed smooth 4-manifold with $\pi_1(M) = \Z^k$ for $k = 1, 2$. Suppose there is an isomorphism\begin{equation}\label{Isomorphism Basic}\xi: H^2(Z)\rightarrow H^2(M)\end{equation}that preserves the cup product and that restricts to a one-to-one correspondence between the basic classes of $Z$ and $M$. Then, the 4-manifold $M$ is irreducible. 
\end{theorem}

\begin{proof}The sole contribution of this proof to the argument due to Park \cite[p. 574]{[Park]} (that uses work of Kotschick \cite{[Kotschick]}) is a slight modification to cover the cases of free abelian fundamental group in the statement of Theorem \ref{Theorem Irreducible Diffeo}. A famous result of Taubes's \cite[Main Theorem]{[Taubes]} says that there is a basic class $k_Z = -c_1(Z, J)\in \mathcal{B}_Z$ given that $Z$ is symplectic, and the hypothesis on the isomorphism (\ref{Isomorphism Basic}) implies that there is at least one basic class $k_M\in \mathcal{B}_M$. Suppose $M = M_1\# M_2$ and, since $\pi_1(M) = \Z^k$, set $\pi_1(M_2) = \Z^k$. Given that the Seiberg-Witten invariant of $M$ satisfies $\Sw_M(k_M)\neq 0$, we have that $b_2^+(M_1)$ and $b_2^+(M_2)$ cannot both be positive \cite[Theorem 2.4.6]{[GompfStipsicz]}. We split the proof in the two cases\begin{enumerate}
\item $b_2^+(M_1) = 0$ and 
\item $b_2^+(M_2) = 0$.
\end{enumerate}
We now consider the first case. If furthermore the second Betti number satisfies $b_2(M_1) = 0$, then $M_1$ is homeomorphic to the 4-sphere by Freedman's theorem \cite[Corollary 1.2.28]{[GompfStipsicz]} and the claim follows. We need to rule out the subcase $b_2(M_1) = b_2^-(M_1) > 0$ and, to do so, we proceed by contradiction. The intersection form over the integers $Q_{M_1}$ of $M_1$ is negative definite. Donaldson's theorem \cite{[Donaldson]} says that $Q_{M_1} \cong n \langle -1 \rangle$ where $n = b_2(M_1) > 0$. Let $\{e_1, \ldots, e_{b_2(M_1)}\}$ be a basis for the second cohomology group $H^2(M_1; \Z)$ with $e_j^2 = -1$ for each $j = 1, \ldots, b_2(M_1)$. The neck pinching argument \cite{[Donaldson], [Kotschick]} implies that there is at least one basic class $k_{M_2}\in \mathcal{B}_{M_2}$ and $\Sw_{M_2}(k_{M_2})\neq 0$. Moreover, the cohomology classes $k_{M_2} + \overset{b_2(M_1)}{\underset{j = 1}{\sum}} a_i e_i$ are basic classes for $M = M_1\cs M_2$, where $a_i = \pm 1$. Consider the basic class $\xi(k_Z)\in \mathcal{B}_{M}$ and set $\xi(k_Z) = k - \overset{b_2(M_1)}{\underset{j = 1}{\sum}} e_j$. As Park observed \cite[p. 574]{[Park]}, the combination $\xi(k_Z) + 2 e_1 = k_{M_2} + e_1 - \overset{b_2(M_1)}{\underset{j = 2}{\sum}} e_j$ is also a basic class of $M$ and $k_Z + 2\xi^{-1}(e_1)\in \mathcal{B}_Z$. A result of Taubes \cite{[Taubes3]} (see \cite[Remark 10.1.16 (b)]{[GompfStipsicz]}) implies that the Poincar\'e dual of $\xi^{-1}(e_1)$ is represented by a symplectically embedded 2-sphere $S\subset Z$ with $[S]^2 = -1$. This is, however, impossible since $Z$ is minimal by hypothesis and we have reached a contradiction. We conclude that $b_2(M_1) = b_2^-(M_1) = 0$ and $M_1$ is homeomorphic to $S^4$. This concludes the proof of the first case.

We now work on the second case $b_2^+(M_2) = 0$. Notice that the previous argument applies verbatim if the intersection form over the integers of $M_2$ is negative definite and, consequently, it rules out the case $b_2^-(M_2) > 0$. We need to rule out the subcase $b_2^-(M_2) = b_2(M_2) = 0$. In this scenario, the Euler characteristic of $M_2$ is zero by our assumption on its fundamental group. We apply the neck pinching argument \cite{[Donaldson], [Kotschick]} once more to conclude that there is at least one basic class $k_{M_1}\in \mathcal{B}_{M_1}$ and $\Sw_{M_1}(k_{M_1})\neq 0$. This, however, is impossible since $\chi(M_2) = 0$ implies that $b_2(M_1)$ is even. The dimension of the moduli space of solutions to the Seiberg-Witten equations and $b_2^+(M_1)$ have opposite parity. In particular, in the case $b_2^+(M_1) = 0 \mod 2$, the Seiberg-Witten invariant of $M_1$ is zero \cite[Section 2]{[Taubes2]}, \cite[Section 2.4]{[GompfStipsicz]}. We conclude that the 4-manifold $M$ is irreducible.

\end{proof}

\section{The second homotopy group of a closed orientable 4-manifold with free abelian fundamental group of small rank.}\label{str}

The main result of this section characterizes the second homotopy group of 4-manifolds with certain free abelian fundamental groups. The main ideas used in the proof of the following result were kindly provided to us by Daniel Kasprowski.

\begin{theorem}\label{Theorem Second Homotopy Group} Let $X$ be a closed oriented 4-manifold with Euler characteristic $\chi = \chi(X)$ and second Betti number $b_2 = b_2(X)$. 
    
    $\bullet$ If the fundamental group of $X$ is infinite cyclic, then its second homotopy group is $\pi_2(X) \cong \Z[\pi_1(X)]^{b_2}$.
    
    $\bullet$ If the fundamental group of $X$ is $\pi_1(X) \cong \Z^2$, then its second homotopy group is $\pi_2(X) \cong \Z \oplus \Z[\pi_1(X)]^{\chi}$. 
    
    In both of these cases, the equivariant intersection form of $X$ is unimodular on the quotient of $\pi_2(X)$ with the radical of the equivariant intersction form (see Definition \ref{defmain}).
\end{theorem}

\begin{proof}We start with the identity in the first clause and assume that $\pi_1(X) \cong \Z$. The universal coefficient spectral theorem \cite[Section 3]{[HKT]} implies that $\pi_2(X)$ is a free $\Z[\pi_1(X)]$-module, whose rank can be computed by using the Cartan-Leray spectral sequence (\cite{[DavisKirk]}). 

We now argue the identity of the second clause. If $\pi_1(X) \cong \Z^2$, the universal coefficient spectral theorem \cite[Section 3]{[HKT]} implies that there is an injective homomorphism \[ i: H^2(\pi_1(X); \Z[\pi_1(X)])\cong \Z \hookrightarrow \pi_2(X)\]
 and \cite[Corollary 4.4]{[HKT]} implies that the quotient of $\pi_2(X)$ by the image of such map is a stably-free $\Z[\pi_1(X)]$-module. Since every stably-free $\Z[\pi_1(X)]$-module is free whenever $\pi_1(X)$ is a finitely generated abelian group (see \cite{[J]}) we have that $\pi_2(X) \cong \Z \oplus \Z[\pi_1(X)]^k$ for some non-negative integer $k$. To conclude that $k=\chi(X)$, we apply a result of Kasprowski-Powell-Teichner \cite[Theorem A]{[KPT]} that says that there are integers $n_1$ and $n_2$ such that there is a diffeomorphism\[ X \cs_{n_1} \CP^2 \cong_{\mathcal{C}^{\infty}} (T^2 \times S^2) \cs_{n_2} \CP^2.\] Since two diffeomorphic 4-manifolds have the same equivariant intersection form and connected summing with a copy of $\CP^2$ increases both the Euler characteristic and the rank of the second homotopy group of a manifold by one, the conclusion follows.

The last part of the statement of Theorem \ref{Theorem Second Homotopy Group} is proven in \cite[Corollary 4.4]{[HKT]}.
\end{proof}

\section{A construction of irreducible 4-manifolds.}\label{Section BK Construction}

\subsection{An irreducible symplectic 4-manifold with free abelian fundamental group}\label{Z}
We describe a construction due to Baldridge-Kirk \cite[Section 4]{[BaldridgeKirk1]}, \cite[Section 4.3]{[BaldridgeKirk]} of a symplectic 4-manifold $M$ (this 4-manifold is denoted by $Z$ in \cite{[BaldridgeKirk]}) that serves as a raw material to produce an infinite family of pairwise non-diffeomorphic irreducible 4-manifolds in the homeomorphism type of $\CP^2 \cs_3 \CPbar^2$; such examples were also constructed by Akhmedov-Park \cite{[AkhmedovPark1]}, Fintushel-Park-Stern \cite{[FintushelParkStern]} and Akhmedov-Park-Baykur \cite{[AkhmedovBaykurPark]}; see Remark \ref{Remark Relation}. Such 4-manifold is  a generalized fiber sum 
\begin{equation}\label{manifold}
M:=M_1 \cs _{\Sigma_2} M_2
\end{equation}
of two symplectic building blocks
 $$M_1:=(T^2 \times T^2) \cs_2 \CPbar^2$$ 
 and $$M_2:=T^2 \times \Sigma_2$$ along a symplectically embedded surface $\Sigma_2$ of genus two with trivial normal bundle. A result of Gompf \cite{[Gompf]} (see \cite[Theorem 10.2.1]{[GompfStipsicz]}) implies that $M$ admits a symplectic structure. Moreover, a result of Usher \cite{[Usher]} allows us to conclude that the symplectic 4-manifold $M$ is minimal and Theorem \ref{Theorem Irreducibility} implies $M$ is irreducible since $\pi_1(M) \cong \Z^6$ is residually finite and $b_2^+(M) > 1$. Its symplectic Kodaira dimension is immediately computed as follows.

 \begin{lemma}\label{Lemma Kodaira Dim}The symplectic Kodaira dimension of the irreducible symplectic 4-manifold $(M, \omega)$ defined in Formula \ref{manifold} is $\Kod(M, \omega) = 2$.
 \end{lemma}

 \begin{proof} Let $K_\omega\in H^2(M)$ be the canonical class of the minimal symplectic 4-manifold $(M, \omega)$. We have that$$K_\omega\cdot K_\omega = 2\chi(M) + 3\sigma(M) > 0$$ and the claim follows from Definition \ref{Definition Kodaira}.
\end{proof}

We now recall some homotopy theoretical traits of the 4-manifold $M$ that are useful for our purposes. Baldridge-Kirk \cite[Theorem 5]{[BaldridgeKirk]} showed that the group $\pi_1(M) \cong \Z^6$ is generated by the based homotopy classes of six loops $x, y, a_1, a_2, b_1, b_2 \subset M$. Moreover, they described  smooth representatives for a basis of the second singular homology group with integer coefficients $H_2(M)\cong \mathbb{Z}^{16}$ in \cite[Proposition 12]{[BaldridgeKirk]}. We refer the reader to their paper for the description of such embedded surfaces, which we will denote by \[H_1, H_2, H_3, F, T_1, R_1, T_2, R_2, T_3, R_3, T_4, R_4, T_1', R_1', T_2', R_2'\] using the same notation. With respect to such basis, the intersection form of $M$ is represented by the matrix

\begin{equation}
Q_M=
    \begin{pmatrix}
-1 & 0 & 0 & 1\\
0 & -1 & 0 & 1\\
0 & 0 & 0 & 1\\
1 & 1 & 1 & 0
\end{pmatrix} 
\bigoplus_{i=1}^6 
\begin{pmatrix}
    0 & 1 \\
    1 & 0
\end{pmatrix}.
\end{equation}
In particular, the six hyperbolic summands correspond to geometrically dual pairs of smoothly embedded $2$-tori. Baldridge and Kirk show that performing torus surgery on the collection of pairwise disjoint 2-tori \begin{equation}\label{tori} \mathcal{T} := \{T_1, T_2, T_3, T_4, T_1', T_2'\} \end{equation} with an appropriate choice of loops and coefficients yields an infinite family of exotic $4$-manifolds with the same homeomorphism type of $\CP^2 \cs_3\CPbar^2$ \cite[Corollary 14]{[BaldridgeKirk]}.

\subsection{Some 4-manifolds obtained from torus surgeries on $M$ and their fundamental groups}\label{2} Baldridge-Kirk perform torus surgeries on the 4-manifold $M$ defined in \ref{manifold} to produce several interesting symplectic 4-manifolds of Euler characteristic six and signature minus two, which include examples with infinite cyclic fundamental group \cite[Theorem 22]{[BaldridgeKirk]} and free abelian fundamental group of rank three in \cite[Corollary 28]{[BaldridgeKirk]}. A modest tweak in their construction yields a myriad of other examples, which include the 4-manifolds of Theorem \ref{Theorem Main} for $k = 3, 5$ \cite{[BaldridgeKirk], [Torres1]}. 

The following result is essentially due to Baldridge-Kirk. Our contribution is to extend their computations and study the fundamental group of the 4-manifolds obtained by torus surgery on the elements of all the possible subcollections of $\mathcal{T}$. In particular, we show that $(-1)$-torus surgery on $M$ on such subfamilies yields a 4-manifold with either a free abelian or a non-abelian fundamdental group.

\begin{theorem}[Baldridge-Kirk \cite{[BaldridgeKirk]}]\label{fundamental groups} Let $M$ be the smooth closed oriented 4-manifold defined in formula \ref{manifold}. The following facts hold.
\begin{enumerate}
    \item Performing $(-1)$-torus surgery on one or two distinct 2-tori of the collection $\mathcal{T}$ produces a closed irreducible symplectic 4-manifold of symplectic Kodaira dimension two and with non-abelian fundamental group.
 
    \item Performing $(-1)$-torus surgery on a collection of three distinct 2-tori $\mathcal{T}_3 \subset  \mathcal{T}$ produces a closed irreducible symplectic 4-manifold $M_{\mathcal{T}_3}$ of symplectic Kodaira dimension two and with fundamental group $\pi_1(M_{\mathcal{T}_3})\cong \Z^3$ provided $\mathcal{T}_3$ is one of the following collections \begin{align*}\{T_1, T_2, T_4\}, \{T_1, T_2, T_1'\}, \{T_1, T_3, T_4\}, \{T_2, T_4, T_2'\}, \\ \{T_2, T_1', T_2'\}, \{T_3, T_4, T_1'\}, \{T_3, T_4, T_2'\}, \{T_3, T_1', T_2'\}.\end{align*} In all the other cases, $\pi_1(M_{\mathcal{T}_3})$ is a non-abelian group.   

    \item Performing $(-1)$-torus surgery on a collection of four distinct 2-tori $\mathcal{T}_4 \subset \mathcal{T}$ produces a closed irreducible symplectic 4-manifold $M_{\mathcal{T}_4}$ of symplectic Kodaira dimension two and with fundamental group $\pi_1(M_{\mathcal{T}_4})\cong \Z^2$ provided $\mathcal{T}_4$ is one of the following collections \begin{align*}\{T_1, T_2, T_3, T_4\}, \{T_1, T_2, T_4, T_1'\}, \{T_1, T_2, T_1', T_2'\}, \\ \{T_1, T_3, T_4, T_1'\}, \{T_1, T_3, T_4, T_2'\}, \{T_1, T_3, T_1', T_2'\}, \\ \{T_2, T_3, T_1', T_2'\}, \{T_2, T_4, T_1', T_2'\}, \{T_3, T_4, T_1', T_2'\}.\end{align*} In all the other cases, $\pi_1(M_{\mathcal{T}_4})$ is a non-abelian group.
    \\The 4-manifold that is obtained if we perform three $(-1)$-torus surgeries and a $p$-torus surgery with $p\neq \pm 1$ on a collection of four distinct 2-tori $\mathcal{T}_4\subset \mathcal{T}$ has fundamental group $\Z^2$.

    \item Performing $(-1)$-torus surgery on a collection of five distinct 2-tori $\mathcal{T}_5 \subset \mathcal{T}$ produces a closed irreducible symplectic 4-manifold of symplectic Kodaira dimension two and with infinite cyclic fundamental group. A generator of such group is given by a loop among $x, y, a_1, a_2,$ $b_1, b_2$ that is contained in $\mathcal{T}\setminus \mathcal{T}_5$.  \\The 4-manifold that is obtained if we perform four $(-1)$-torus surgeries and a $p$-torus surgery with $p\neq \pm 1$ on a collection of five distinct 2-tori $\mathcal{T}_5\subset \mathcal{T}$ has infinite cyclic fundamental group.
    
\end{enumerate}
    
\end{theorem}

\begin{proof} The claims regarding irreducibility and the existence of a symplectic structure of symplectic Kodaira dimension two follow from Theorem \ref{Theorem Luttinger Surgery}, Theorem \ref{Theorem Irreducibility} and Lemma \ref{Lemma Kodaira Dim}. The proof of the other claims is based heavily on Baldridge-Kirk's computation of the fundamental group of of the complement $C$ of the disjoint union of the $2$-tori in the collection $\mathcal{T}$ in \cite[Theorem 11]{[BaldridgeKirk]} and of the based homotopy classes of their meridians, together with Seifert-van Kampen's theorem. In particular, \cite[Theorem 11]{[BaldridgeKirk]} tells us that the $\pi_1(C)$ is generated by the based homotopy classes of the loops $x, y, a_1, a_2, b_1, b_2$ with relations
   \begin{equation}\label{rel comp}
   \begin{split}
       &[x,a_1]=1, [x, a_2]=1, [y,a_1]=1, [y,a_2]=1, [x,y]=1
       \\ &[b_1,b_2]=1, [a_1,b_1]=1, [a_1,b_2]=1, [a_2,b_2]=1.
    \end{split}
   \end{equation}

 The rows of Table 1 contain the 2-tori of the collection $\mathcal{T}$ together with the based homotopy classes of their meridians, the surgery curves and the relations imposed on the generators of $\pi_1(C)$ by gluing back the 2-tori while performing a $p$-torus surgery. 
\begin{center}
\begin{table}
\label{table1}
   \begin{tabular}{|c|c|c|c|}
\hline
2-torus & Meridian & Surgery curve & $p$-torus surgery relation\\
\hline
$T_1$ & $[b_1^{-1},y^{-1}]$ & $x$ & $x=[b_1^{-1},y^{-1}]^{-p} $ \\
\hline
$T_2$ & $[x^{-1},b_1]$ & $a_1$ & $a_1=[x^{-1},b_1]^{-p}$ \\
\hline
$T_3$ & $[b_2^{-1},y^{-1}]$ & $a_2$ & $a_2=[b_2^{-1},y^{-1}]^{-p}$ \\
\hline
$T_4$ & $[x^{-1},b_2]$ & $y$ & $y=[x^{-1},b_2]^{-p}$ \\
\hline
$T_1'$ & $[a_2^{-1},a_1^{-1}]$ & $b_1$ & $b_1=[a_2^{-1},a_1^{-1}]^{-p}$   \\
\hline
$T_2'$ & $[b_1,a_2]$ & $b_2$ & $b_2=[b_1,a_2]^{-p}$ \\
\hline
\end{tabular}
\vspace{3mm}
\caption{2-tori, meridians, surgery curves and relations introduced to $\pi_1$.}
\end{table}
\end{center}

One can compute the free abelian fundamental groups in the statement by following the strategy employed in \cite[Theorem 13]{[BaldridgeKirk]}. We outline the computations for the case $\mathcal{T}_4=\{T_3, T_4, T_1', T_2'\}$ for the convenience of the reader. The fundamental group of the 4-manifold $M_{T_3, T_4, T_1', T_2'}$ is generated by the loops $x, y, a_1, a_2, b_1, b_2$ with relations \ref{rel comp} and
\begin{equation}\label{r}
\begin{split}
& a_2=[b_2^{-1}, y^{-1}], y=[x^{-1}, b_2], b_1=[a_2^{-1},a_1^{-1}], b_2=[b_1,a_2],\\ 
&[b_1,y]=1, [b_1,x]=1. 
\end{split}
\end{equation}

Since $x$ commutes both with $b_1$ and $a_2$ by \ref{rel comp} and \ref{r}, one gets that $y=1$ by substituting the expression for $b_2$ into the one for $y$ in \ref{r}. This in turn implies that $a_2=1, b_1=1$ and $b_2=1$, while $x$ and $a_1$ commute by \ref{rel comp}. It hence follows that $$\pi_1(M_{T_3, T_4, T_1', T_2'}) \cong \langle x, a_1 \mid [x,a_1]=1 \rangle \cong \Z^2.$$

The non-abelian cases are carried out by first finding a presentation of the fundamental group of the 4-manifold obtained by torus surgery on $M$ by means of generators and relations using Seifert-van Kampen's theorem and then building a surjective group homomorphism onto the quaternion group $$Q_8:=\{ \pm 1, \pm i, \pm j, \pm k\}\subset \mathbb{H}.$$ 

We now report the computations of the non-abelian fundamental group of the manifold $M_{T_1,T_2,T_3,T_1'}$ as a leading example. Using the Seifert-van Kampen's theorem we conclude that the fundamental group $\pi_1(M_{T_1,T_2,T_3,T_1'})$ is generated by the based homotopy classes of the loops $x, y, a_1, a_2,$ $ b_1, b_2$. The relations are the ones in (\ref{rel comp}), characterizing the fundamental group of the complement of the union 2-tori, plus the following ones:

\begin{equation}\label{rel}
    \begin{split}
        & x=[b_1^{-1},y^{-1}], a_1=[x^{-1}, b_1], a_2=[b_2^{-1}, y^{-1}]\\ 
        &b_1=[a_2^{-1},a_1^{-1}], [b_1,a_2]=1, [x,b_2]=1.
    \end{split}
\end{equation}
Substituting the expression for $b_1$ in (\ref{rel}) into the one for $a_1$ and using the fact that $x$ commutes both with $a_1$ and $a_2$ by (\ref{rel comp}), we get that $a_1=1$. By substituting $a_1$ into the relations (\ref{rel}), we get that also $b_1=1$ and hence $x=1$. In particular, this implies that the fundamental group $\pi_1(M_{T_1,T_2,T_3,T_1'})$ is generated by the based homotopy classes of $y$ and $b_2$. The non-trivial relations are
\begin{equation}
    [y,[b_2^{-1},y^{-1}]]=1, [[b_2^{-1},y^{-1}], b_2]=1
\end{equation}
and one can check that they are equivalent to the conditions
\begin{equation}
    [y,b_2]=[b_2,y^{-1}]=[b_2^{-1},y].
\end{equation}
In particular, we have just proved that $$\pi_1(M_{T_1,T_2,T_3,T_1'})\cong \langle y, b_2 \mid [y,b_2]=[b_2, y^{-1}]=[b_2^{-1},y] \rangle.$$ This group is non-abelian, since it is possible to define a group homomorphism 
\begin{equation}
\phi: \pi_1(M_{T_1,T_2,T_3,T_1'}) \to Q_8
\end{equation}
by setting $\phi(y)=i$ and $\phi(b_2)=j$. 
\end{proof}

\begin{remark}\label{Remark Relation} Exotic irreducible smooth structures on $\CP^2\cs_3\CPbar^2$ were constructed by Akhmedov-Park \cite{[AkhmedovPark1]} and Baldridge-Kirk \cite{[BaldridgeKirk1], [BaldridgeKirk]} by performing six torus surgeries to the 4-manifold described in formula (\ref{manifold}). Akhmedov-Baykur-Park \cite{[AkhmedovBaykurPark]} showed that the difference in these two constructions is the choice of three out of the six 2-tori involved in the surgeries. In particular, the mechanism described in Section \ref{SW} applies verbatim to the 4-manifolds of Theorem \ref{fundamental groups}. 
\end{remark}

In each item of Theorem \ref{fundamental groups} there are various choices for the sub-collection of tori to perform the surgery.
\begin{question}
    Are the resulting manifolds diffeomorphic to each other?
\end{question}

\section{A surgical procedure to represent homology classes by surfaces with trivial $\pi_1$-inclusion.}\label{construction}



We work in the following setting throughout this section.
\begin{itemize}
    \item The 4-manifold $X'=X_{T,\gamma}(p)$ is the result of a $p$-torus surgery to a closed smooth oriented 4-manifold $X$ on a homologically essential $2$-torus $T\subset X$ along a primitive curve $\gamma \subset T$ and $T' \subset X'$ is the core 2-torus after the surgery.
    \item There is a smoothly embedded 2-torus $S \subset X$ that is nowhere tangent to $T$ satisfying $$\gamma=S \cap T=\sigma_1,$$ where $\sigma_1, \sigma_2 \subset S$ are two simple loops that generate $\pi_1(S)$.
    \item For a given framing $\nu(T) \cong \gamma \times \eta \times D^2$, the intersection $S \cap \nu(T)$ is of the form $\gamma \times \{p_{\eta}\} \times D^1$, where $D^1 \subset D^2$ is the standard interval;
    \item There is a smoothly embedded closed oriented genus two surface $\Sigma\subset X$ with trivial normal bundle and that is geometrically dual to the 2-torus $T$. Moreover, $[\sigma_2] \notin (i_{\Sigma})_*(\pi_1(\Sigma))$, where $i_{\Sigma}: \Sigma \hookrightarrow X$ is the inclusion map.
\end{itemize}

  \begin{figure}
  \includegraphics[width= 0.8\textwidth, trim=170 2000 170 0, clip]{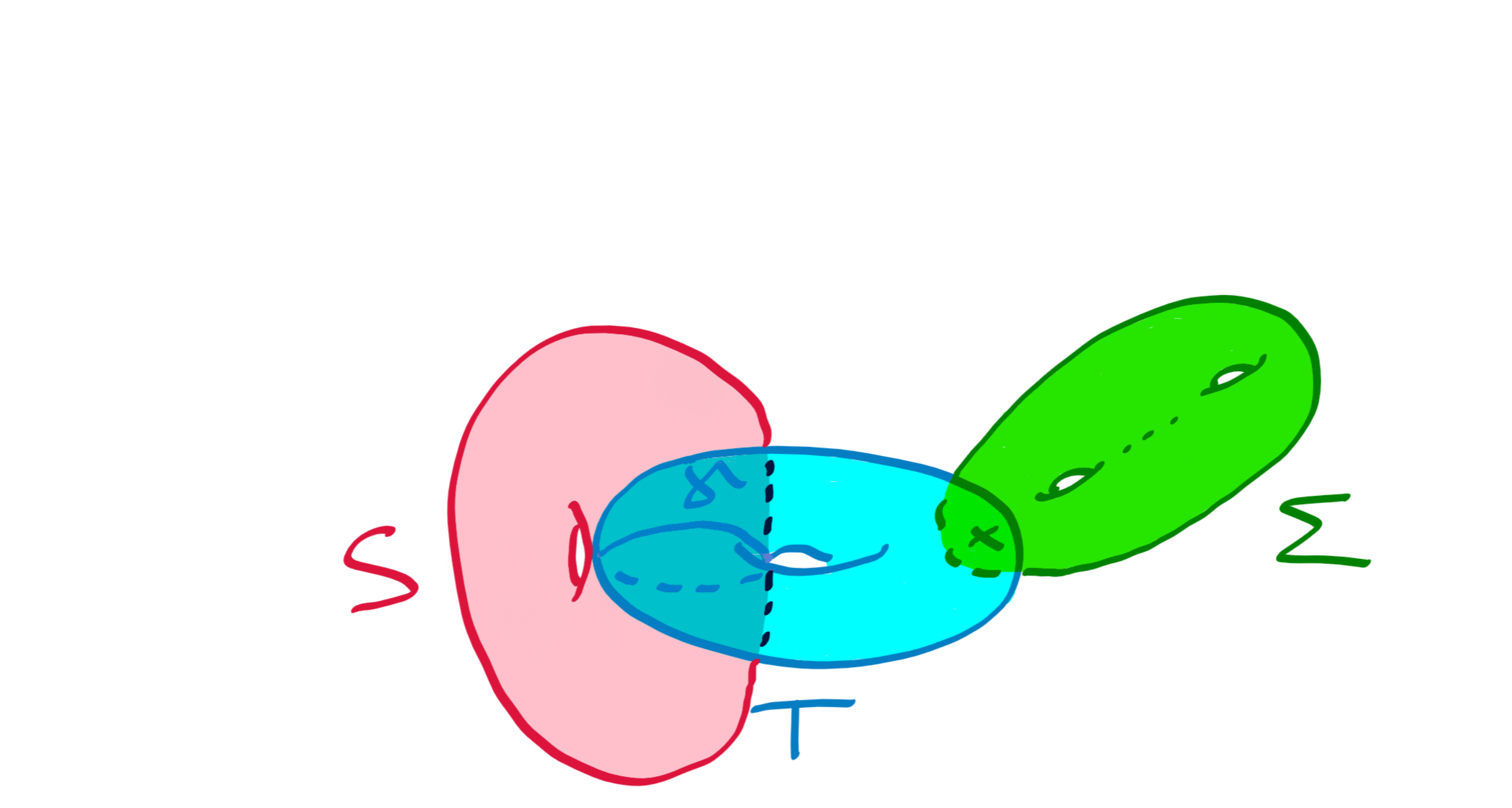} 
  \caption{Reciprocal positions of $T$, $S$ and $\Sigma$ under the hypothesis of Section \ref{construction}.}
  \end{figure}

 Notice that the assumption on the intersection $S \cap \nu(T)$ implies that, up to an isotopy inside $X$, we can make $S$ disjoint from the surgery 2-torus $T$, so that it defines a smoothly embedded surface $$S' \subset X \setminus \nu(T) =  X' \setminus \nu(T') \subset X'$$  with curves $\sigma_1', \sigma_2' \subset S'$ corresponding to $\sigma_1, \sigma_2 \subset S$. In particular, since $S \cap \nu(T)= \gamma \times \{p_{\eta}\} \times D^1$, one can do so by moving the standard interval $D^1 \subset D^2$ along one normal direction, until it degenerates to one point and then disappears.
 
 \begin{lemma}\label{lemma}
    Under the above assumptions, there is a smoothly embedded closed connected oriented surface $S^*\subset X'$ of genus $2pg$ such that $[S^*]=[S'] \in H_2(X')$ and $[\sigma_2'] \notin i_*(\pi_1(S^*))\subset \pi_1(X')$, where $i_*: \pi_1(S^*) \to \pi_1(X')$ is the homomorphism induced by the inclusion. 
\end{lemma}

\begin{proof}
        
    Consider a properly embedded cylinder $C:= S \setminus \nu(T) \subset  X\setminus \nu(T)$ that bounds two parallel push-offs of $\gamma$. Up to a slight perturbation, we can suppose that such push-offs are of the form $$\gamma^i=\gamma \times \{p_{\eta,i}\} \times \{p_{\mu, i}\} \subset \gamma \times \eta \times \partial D^2$$ inside $\partial \nu(T)\cong \gamma \times \eta \times \partial D^2$, where $p_{\eta,1} \neq p_{\eta,2}$.  Since $X\setminus \nu(T)=X' \setminus \nu(T')$, such cylinder is also contained in $X'$, where it is denoted by $C' \subset X'$. Let $\{\Sigma_{j}^1, \Sigma_{j}^2 \mid j=1, \dots, p\}$ be $2p$ distinct parallel copies of the properly embedded surface $$\Sigma \setminus \nu(T)\subset X\setminus \nu(T).$$ In particular, we can suppose without loss of generality that 
    $$\partial \Sigma_{j}^i=\{ p_{\gamma,j}\} \times \{p_{\eta,i}\} \times \partial D^2 \subset \gamma \times \eta \times \partial D^2$$ for $i=1,2, j=1, \dots, p$, where $p_{\gamma,h} \neq p_{\gamma,k}$ for any $h \neq k$. The union of the cylinder $C$ with the surfaces $\Sigma_{j}^i$ gives a properly immersed surface $$\tilde S:=C \cup_{i,j} \Sigma^i_j \subset X \setminus \nu(T)$$ with $\partial \tilde S \subset T_1 \sqcup T_2$, where for $i=1,2$ we set $$T_i:=\gamma \times \{p_{\eta,i}\} \times \partial D^2\subset \partial \nu(T).$$ In particular, $\partial \tilde S$ is the union of the push-offs $\gamma^i$ and $\partial \Sigma^i_{j}$ for $j=1, \dots, p$, which inside the 2-tori $T_1, T_2$ respectively look like the union of a longitude with $p$ disjoint parallel copies of a meridian curve, given by $\partial \Sigma^i_j$ for $j=1, \dots,p$, as in Figure \ref{curves/fig}. We orient $\gamma_i$ as a boundary component of $C \subset S$, while an orientation on $\Sigma^i_j$ is defined by imposing that $\gamma_i$ and $\partial \Sigma^i_j$ have both the same or the opposite orientation as the curves $\gamma$ and $\mu_T$ inside $T_i$. Notice that in this way we have that $[\Sigma^1_j]=-[\Sigma^2_k]\in H_2(X,\nu(T))$ for all $j,k=1, \dots,p$.

     \begin{figure}
  \includegraphics[width= 0.8\textwidth, trim=170 2000 170 0, clip]{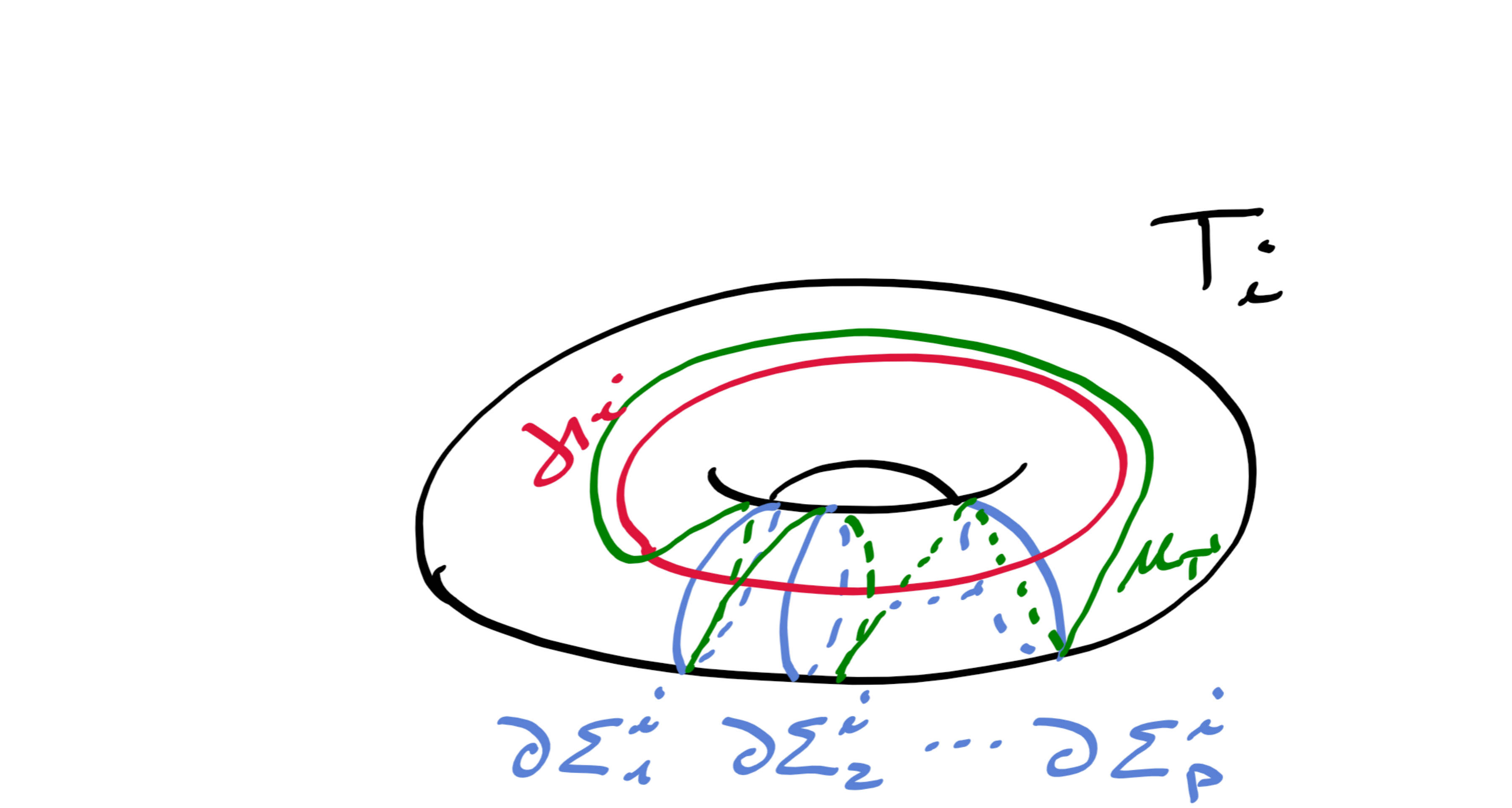} 
  \caption{Curves in the 2-tori $T_1$ and $T_2$.}
  \label{curves/fig}
  \end{figure}
    
    There are hence $2p$ points of self-intersection between the cylinder $C$ and the surfaces $\Sigma_{j}^i$ inside the two 2-tori $T_1$ and $T_2$. We can resolve such singularities by gluing bands in the following way. Locally around the intersection point between $C$ and $\Sigma^i_j$, we picture the two surfaces in a 3-dimensional chart with coordinates $x,y,z$ in such a way that $C$ is locally sent to the the half-plane $H_1:=\{x=0, y \geq 0\}$, $\Sigma^i_j$ goes onto $H_2:=\{z=0,y \leq 0\}$ and the intersection point $C \cap \Sigma^i_j$ is sent to the origin. Apply a small translation to $H_1$ along the y-axis to make it disjoint from $H_2$, e.g. by mapping it to the subset $H_1':=\{x=0, y \geq \epsilon\}$ with $\epsilon >0$, and glue $H_1'$ and $H_2$ by means of a 2-dimensional band doing a $\frac{\pi}{2}$-twist in the $z$-direction, where the sense of rotation is determined by respecting the local orientations of $C$ and $\Sigma^i_j$. Up to composing with an ambient isotopy of $\mathbb{R}^3$, we now get a surface with boundary contained in the plane of equation ${y=0}$ as in Figure \ref{Figure4}.
    
     \begin{figure}
  \includegraphics[width= 0.8\textwidth, trim=170 1900 170 0, clip]{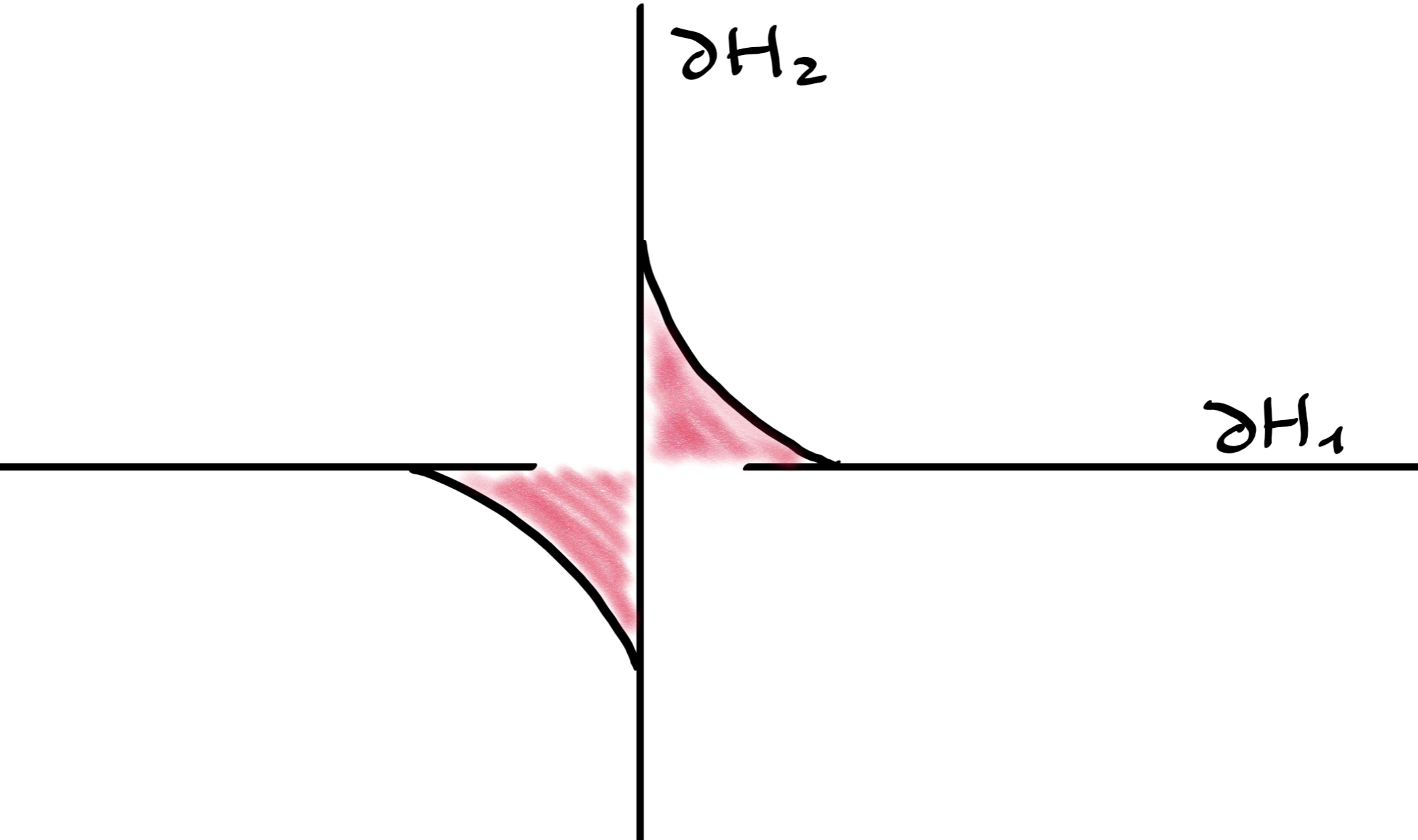} 
  \caption{Band desingularization.}
  \label{Figure4}
  \end{figure}
    
    By locally substituting the intersection point of $C$ with $\Sigma^i_j$ with this local model for every $j=1, \dots, p$, we get a new properly embedded surface $$\widehat{S} \subset X \setminus \nu(T)=X'\setminus \nu(T')$$ of genus $2gp$ with two boundary components which are isotopic inside each $T_i$ to an oriented closed curve which is the connected sum of $\gamma^i$ with the $\partial \Sigma^i_j$'s. In particular, such closed curves are isotopic to the meridian of the core 2-torus $T'$ and have opposite orientations. This allows us to cap-off the boundary components of $\widehat{S} \subset X'$ with two 2-disks bounding two parallel copies of the meridian of $T'$ inside $\nu(T')$ and finally get the desired submanifold $S^*\subset X'$.

     Let us now consider the homology class $[S^*] \in H_2(X')$. Since the core 2-torus $T'\subset X'$ is null-homologous by Proposition \ref{surgery}, the natural map $H_2(X') \to H_2(X',\nu(T'))$ is injective and $[S^*]=[S']\in H_2(X')$ if and only if $[S^*]=[S']\in H_2(X',\nu(T'))$. This last identity holds since $$[S^*]=[\widehat{S} \cup D^2_1 \cup -D^2_2]=[\widehat{S}]=[C' \cup_{i,j} \Sigma^i_j]=[C']=[S']\in H_2(X',\nu(T'))$$ given that we have oriented each surface $\Sigma^i_j$ in such a way that $[\Sigma^1_j]=-[\Sigma^2_j]\in H_2(X',\nu(T'))$ for every $j=1, \dots, p$, while $D^2_1, D^2_2$ are the two disks used to cap-off the boundary components of $\widehat{S}$.
\end{proof}

\section{Topological prototypes.}\label{Section Infinite Families}
We now pin down the homeomorphism type of the 4-manifolds that were constructed in Section \ref{2} and begin with those with infinite cyclic fundamental group. For 2-tori $T \in \mathcal{T}:= \{T_1, T_2, T_3, T_4, T_1', T_2'\}$ and $p \in \N$, denote by $M_{T}(p)$ the $4$-manifold that is obtained from $M$ defined in formula (\ref{manifold}) by performing $(-1)$-torus surgery on four distinct elements of the set $\mathcal{T} \setminus \{T\}$ and a $p$-torus surgery on the remaining $2$-torus, with the choices surgery curves explained at the beginning of Section \ref{2}. Theorem \ref{fundamental groups} says that $\pi_1(M_{T}(p)) \cong \Z$ for every $p \in \N$. As it was mentioned in the introduction, these 4-manifolds need not be new \cite{[BaldridgeKirk], [Torres1], [Torres2]}.
    
\begin{theorem}\label{main}    
    For every $p \in \N$, there exists a homeomorphism $$M_{T}(p) \cong_{\mathcal{C}^0} \cs_2\CP^2 \cs_4 \CPbar^2 \cs (S^1 \times S^3).$$
\end{theorem}
\begin{proof} For the sake of clarity in the exposition, we argue in detail the case where $T = T_2$ and the $p$-torus surgery is performed along the $2$-torus $T_4$. The proof for all other choices is carried out in a similar way. Theorem \ref{fundamental groups} implies that $$\pi_1(M_{T_2}(p)) = \langle a_1 \rangle \cong \mathbb{Z}.$$ In order to apply Theorem \ref{classification} to pin down the homeomorphism type of $M_{T_2}(p)$, we compute its equivariant intersection form up to isometries.

    We first consider the intersection form of $M_{T_2}(p)$ over the integers. Proposition \ref{surgery} implies that $b_2(M_{T_2}(p))=b_2(M)-2 \cdot 5=6$. In particular, the second homology group $$H_2(M_{T_2}(p))\cong \mathbb{Z}^6$$ is generated by the surfaces $T_2, R_2, H_1,$ $H_2, H_3, F$ (see \cite[Proposition 12]{[BaldridgeKirk]}). Notice that these surfaces do not intersect the 2-tori in $M$ that are used in the surgeries to produce $M_{T_2}(p)$ and, hence, they define six embedded surfaces inside $M_{T_2}(p)$ for every $p \in \N$. It follows that the intersection form of $M_{T_2}(p)$ over the integers coincides with the one of $\cs_2 \CP^2 \cs_4 \CPbar^2 \cs (S^1 \times S^3)$ for every $p \in \N$. To finish the proof of Theorem \ref{main}, we show that the equivariant intersection form of $M_{T_2}(p)$ is extended from the integers.

    The surfaces $H_1, H_2, H_3$ and $R_2$ admit a lift to the universal cover of $M_{T_2}(p)$, since the image of their fundamental groups via the inclusion induced homomorphism is trivial. Moreover, for a fixed choice of a lifting, these surfaces define elements of the second homology group of the universal cover of $M_{T_2}(p)$, on which the equivariant intersection form can be computed. In the following, we will build two new surfaces that are homologous to $T_2$ and $F$, respectively, and which also admitting a lift to the universal cover of $M_{T_2}(p)$. We start with $T_2$. Apply Lemma \ref{lemma} with the choices $S=T_2$, $\Sigma=R_4$, $T=T_4$, $\gamma=y$, $X'=M_{T_2}(p)$ and $X$ equal to the 4-manifold obtained by performing torus surgery on the 2-tori $T_1, T_3, T_1', T_2'$ to obtain a smoothly embedded oriented surface $T_2^{*p}\subset M_{T_2}(p)$ such that $i_*(\pi_1(T_2^{*p}))=0 \subset \pi_1(M_{T_2}(p))$, where $i_*$ is the inclusion-induced homomorphism. In order to find a new representative for $[F] \in H_2(M_{T_2}(p))$, we consider the connected sum decomposition of $F$ into two 2-tori $F_1$ and $F_2$ as in \cite[Section 4.3]{[BaldridgeKirk]}, where $\pi_1(F_i) \cong \langle a_i, b_i \mid [a_i,b_i]=1 \rangle \cong \Z^2$ for $i=1,2$. We will replace $F_1$ in such decomposition with a surface that admits a lift to the universal cover of the ambient manifold. By the definition of the map used in (\ref{manifold}) for gluing the two building blocks $M_1$ and $M_2$, we have that (up to an ambient isotopy of $M_{T_2}(p)$) $F_1$ intersects the tubular neighbourhood $\nu(T_1')$ along a cylinder bounding two parallel push-offs of the loop $b_1$, which is killed at the level of the fundamental group by the torus surgery we perform on $T_1'$. Notice that one might first think to set $S=F_1$, $\Sigma=R_1'$, $T=T'_1$, and $\gamma=b_1$ in order to apply Lemma \ref{lemma}, and get a new liftable surface representing the same homology class as $F_1$. However, in this case the procedure described in Lemma \ref{lemma} would build a surface that cannot be lifted to the universal cover of $M_{T_2}(p)$, since the inclusion $R_1' \subset M_{T_2}(p)$ does not induce the trivial homomorphism between fundamental groups.  We can get around this issue by iterating the process and applying Lemma \ref{lemma} with a different choice of $\Sigma$. In particular, for our purposes we can set $\Sigma$  to be the geometrically dual surface to $T_1'$ obtained by applying Lemma \ref{lemma} to $R_1'$ with choices $S=R_1'$, $\Sigma=R_3$ and $T=T_3$. Indeed, a basis of $\pi_1(R_1')$ is given by $a_1, a_2$ and $a_2$ is null-homotopic in $M_{T_2}(p)$ thanks to the torus surgery on $T_3$. Moreover, $T_3$ has a geometrically dual 2-torus $R_3$ with the desired properties. In this way, we get out of $F_1$ a new surface $F_1^{*p}$ in the same homology class. The surface $F^{*p}$ obtained by substituting $F_1$ with $F_1^{*p}$ in the connected sum decomposition of $F$ gives the desired new representative of $[F]$.
    
    We conclude that the equivariant intersection form of $M_{T_2}(p)$ can be computed by using homology classes in $H_2(\widetilde{M}_{T_2}(p)) \cong \pi_2(M_{T_2}(p))$ that are represented by lifts of the surfaces $H_1$, $H_2$, $H_3$, $F^{*p}$, $T_2^{*p}$ and $R_2$. In particular, we point out that there is a handlebody decomposition of $M_{T_2}(p)$ for which the hypothesis of Example \ref{example} are satisfied and hence the equivariant intersection form if $M_{T_2}(p)$ is extended from the integers. Indeed, starting from the fiber sum decomposition (\ref{manifold}) of $M = M_1 \cs _{\Sigma_2} M_2$, one can also see $M_{T_2}(p)$ as a generalized fiber sum 
\begin{equation*}
    M_{T_2}(p)=M_{1,T_2}(p) \cs _{\Sigma_2} M_{2,T_2}(p) 
\end{equation*}along a surface of genus two, where $M_{1,T_2}(p)$ is obtained from $M_1$ by $(-1)$-torus surgery on $T_1'$ and $T_2'$ and $M_{2,T_2}(p)$ is defined by simultaneously performing two $(-1)$-torus surgeries on $T_1$ and $T_3$ and a $p$-torus surgery on $T_4$ to $M_2$. The above decomposition implies that $M_{T_2}(p)$ is obtained by performing four loop surgeries and a 2-sphere surgery to the connected sum $M_{1,T_2}(p) \cs M_{2,T_2}(p)$. Moreover, all such surgery operations happen away from the six embedded surfaces we want to use to compute the equivariant intersection form. We can hence build a handlebody decomposition of $M_{T_2}(p)$ starting from the standard ones of $M_1=T^4 \cs_2 \CPbar^2$ and $M_2=T^2 \times \Sigma_2$ (see \cite[Chapter 4]{[Akbulut]}), by keeping track of the result of the loop surgeries, the sphere surgery and all the torus surgery needed to construct $M_{T_2}(p)$. In particular, the single 1-handle giving a non-trivial element of $\pi_1(M_{T_2}(p)) \cong \langle a_1 \rangle$ will come from the 1-handle in the handlebody decomposition of $T^2 \times \Sigma_2$ corresponding to the loop $a_1$ in the symplectic basis $\{a_1, b_1, a_2, b_2\}$ of $\pi_1(\Sigma_2)$. It is then not hard to verify that $H_1, H_2, H_3, F^{*p}, T_2^{*p}$ and $R_2$ are contained in the complement of such 1-handle inside the 2-skeleton of $M_{T_2}(p)$ and the conclusion follows.
\end{proof}
\begin{remark} If in the construction in the proof of Theorem \ref{main}, we perform a p-torus surgery with $p\neq 1$ instead of $p = 1$, then the modified representatives of the homology classes under consideration will have larger genus. However, their algebraic intersection will not change. 
\end{remark}

We spend the remaining part of this section working out the homeomorphism class of the 4-manifolds with free abelian group of rank two that are under consideration. Let us now pick a pair of distinct 2-tori$$T, T' \in \mathcal{T}=\{T_1, T_2, T_3, T_4, T_1', T_2'\}$$ and let $M_{T, T'}(p)$ be the 4-manifold obtained by performing $(-1)$-torus surgery on three distinct 2-tori in $\mathcal{T} \setminus \{T, T'\}$ and a $p$-torus surgery on the remaining 2-torus to the 4-manifold constructed in formula (\ref{manifold}) (cf. \cite{[AkhmedovPark1], [BaldridgeKirk], [Torres1], [Torres2]}). Theorem  \ref{fundamental groups} implies that we can make these choices in such a way that $\pi_1(M_{T, T'}(p)) \cong \Z^2$.

\begin{theorem}\label{main2}    
    For every $p \in \N$, there exists a homeomorphism $$M_{T,T'}(p) \cong_{\mathcal{C}^0} \cs_2\CP^2 \cs_4 \CPbar^2 \cs (T^2 \times S^2).$$
\end{theorem}

\begin{proof}The strategy of the proof is the same as the one of  Theorem \ref{main}: we will show that the equivariant intersection form of $M_{T,T'}(p)$ is extended from the integers. For the sake of clarity, we consider the case with choices $T = T_1$ and $T' = T_2$, and perform a $p$-torus surgery on $T_4$. The other cases follow from a similar argument. Theorem \ref{fundamental groups} implies that the group $\pi_1(M_{T_1,T_2}(p))\cong \Z^2$ is generated by the based homotopy classes of $x$ and $a_1$. On the other hand, as already done in the proof of Theorem \ref{main}, one can find new representatives for the second homology classes given by the six surfaces $H_1, H_2, H_3, F, T_2, R_2$ that do lift to the universal cover. On such elements, the equivariant intersection form is given by the matrix  
    \begin{equation}
Q_p=
    \begin{pmatrix}
-1 & 0 & 0 & 1\\
0 & -1 & 0 & 1\\
0 & 0 & 0 & 1\\
1 & 1 & 1 & 0
\end{pmatrix} 
\bigoplus 
\begin{pmatrix}
    0 & 1 \\
    1 & 0
\end{pmatrix}.
\end{equation}

Since $H_2(\widetilde{M}_{T_1,T_2}(p))/R(\widetilde{I}_{M_{T_1,T_2}(p)})$ is a free $\Z[\pi_1(M_{T_1,T_2}(p))]$-module of rank $\chi(M_{T_1,T_2}(p))=6$ by Theorem, \ref{Theorem Second Homotopy Group}, the fact that $Q_p$ is unimodular together with Remark \ref{reduced} implies that $$\widetilde I_{M_{T_1,T_2}(p)} \cong Q_p \otimes_{\Z} \Z[\pi_1(M_{T_1,T_2}(p))]$$ 
for every $p \in \N$. Since both $M_{T_1,T_2}(p)$ and its universal cover are non-spin for any $p \in \N$, Theorem \ref{main2} follows from Theorem \ref{classification}.
\end{proof}

\section{Stabilization of $p$-torus surgeries.}\label{stab}

We extend results of Baykur-Sunukjian \cite{[BaykurSunukjian]} in the following proposition, which will be used in the proof of Item (3) in Theorem \ref{Theorem Main}. Our extension is based on coupling Baykur-Sunukjian's work with a well-known trick due to Moishezon's \cite{[Moishezon]} that describes a relation between 0-torus surgeries and surgery along loops. 

 \begin{proposition}\label{Proposition S} Let $X$ be a closed smooth oriented 4-manifold that contains a smoothly embedded 2-torus $T \subset X$ with trivial normal bundle and let $X_{T, \gamma}(p)$ be the 4-manifold defined in (\ref{Torus Surgery}), where the surgery curve $\gamma \subset T$ is a generator of the fundamental group $\pi_1(T)$. Suppose that $X\setminus \nu(T)$ is non-spin and that the inclusion $X\setminus \nu(T)\rightarrow X$ induces an isomorphism of fundamental groups $\pi_1(X\setminus \nu(T))\rightarrow \pi_1(X)$. There is a diffeomorphism \begin{center}$X_{T, \gamma}(p) \cs {(S^2 \times S^2)}\cong_{\mathcal{C}^{\infty}} X_{T, \gamma}(0) \cs {(S^2 \times S^2)}$\end{center} for every $p \in \N$. 
 \end{proposition}

 \begin{proof} There is a cobordism\begin{equation}\label{Cobordism 1}V = (I\times X) \cup R_2\end{equation}for a five-dimensional round 2-handle $R_2 = S^1\times D^2\times D^2$ such that $\partial_-V = X$ and $\partial_+ V = X_{T, \gamma}(p)$ by a result of Baykur-Sunukjian \cite[Lemma 1]{[BaykurSunukjian]}. Asimov's fundamental lemma of round handles \cite{[Asimov]} says that the round 2-handle $R_2$ decomposes into a 2-handle $h_2 = D^2\times D^3$ and a 3-handle $h_2 = D^3\times D^2$ that are attached independently, and the cobordism (\ref{Cobordism 1}) can be expressed as
 \begin{equation}\label{Cobordism 2}V = (I\times X) \cup h_2 \cup h_3.\end{equation} We now look in further detail at the middle level $V_{1/2}$ of (\ref{Cobordism 2}) in order to prove Proposition \ref{Proposition S} as it was done by Baykur-Sunukjian in \cite{[BaykurSunukjian]}. The attaching circle of the five-dimensional 2-handle in (\ref{Cobordism 2}) is the surgery curve $\gamma\subset X$ (cf. \cite[Proof of Lemma 4]{[BaykurSunukjian]}). In particular, we have\begin{equation}\label{Loop Surgery 1}V_{1/2} = (X\setminus \nu(\gamma))\cup (D^2\times S^2),\end{equation}where the loop $\gamma\subset T\subset X$ is framed using the product framing of the 2-torus $T$. Apply Moisheson's trick \cite{[Moishezon]} (see \cite{[Gompf2]}) to describe (\ref{Loop Surgery 1}) as the result of doing surgery along a null-homotopic loop $\gamma'\subset X_{T, \gamma}(0)$ and conclude that\begin{equation}V_{1/2} = X_{T, \gamma}(0)\cs {(S^2\times S^2)}\end{equation} by using the hypothesis that the complement $X\setminus \nu(T)$ is non-spin. Notice that there is a diffeomorphism \[X_{T, \gamma}(0)\cs {(S^2\times S^2)} \cong_{\mathcal C^\infty} X_{T, \gamma}(0)\cs {(S^2\simtimes S^2}),\] where $S^2\simtimes S^2$ is the twisted $S^2$-bundle over $S^2$.

 Turn the cobordism (\ref{Cobordism 2}) upside down, to see that $V_{1/2}$ is also obtained by attaching a five-dimensional 2-handle to $X_{T, \gamma}(p)$ along a null-homotopic loop: this loop is null-homotopic since we have assumed that the meridian $\mu_T$ of the 2-torus $T\subset X$ is null-homotopic in $X\setminus \nu(T)$. Given that $X_{T, \gamma}(p)\setminus \nu(\widehat{T}_{p})=X \setminus \nu(T)$ is non-spin, we conclude that\begin{equation}V_{1/2} = X_{T, \gamma}(p)\cs {(S^2\times S^2)}\end{equation} and the proposition follows.

 \end{proof}

\section{Proof of Theorem \ref{Theorem Main}.}\label{Section Proof}

We provide a detailed proof of the statement in the case $k = 3$. Set $M_{p,3}\coloneq M_T(p)$ and $N_{p,3} \coloneq M_{T,T'}(p)$, where $M_T(p)$ is the 4-manifold defined in the proof of Theorem \ref{main}, and $M_{T, T'}(p)$ is the 4-manifold defined the proof of Theorem \ref{main2}; cf. Theorem \ref{fundamental groups}. Theorem \ref{Theorem Smooth Structures} along with the discussion in Section \ref{SW} and Theorem \ref{Theorem Irreducible Diffeo} imply that the two collections\begin{center}$\{M_{p,3} : p \in \N\}$ and $\{N_{p,3} : p \in \N\}$\end{center}consist of pairwise non-diffeomorphic elements that are irreducible provided $p\neq 0$; see Remark \ref{Remark Relation}. Their homeomorphism classes have been determined in Theorem \ref{main} and Theorem \ref{main2}. The third clause of Theorem \ref{Theorem Main} follows from Proposition \ref{Proposition S}. The 4-manifolds $M_{1,3}$ and $N_{1,3}$ admit a symplectic structure of symplectic Kodaira dimension two by \ref{Theorem Luttinger Surgery Kodaira}, in view of \ref{Lemma Kodaira Dim}.

The infinite collections for the values $k\in \{2, 4, 5\}$ are readily available in the literature \cite{[AkhmedovPark2], [BaldridgeKirk], [Torres1]} as we have summarized in Table 2. From left to right, the columns of such table contain the value of $k \in \{ 2,3,4,5\}$, the Euler characteristic and the signature of the raw material in the fourth column. Starting from such 4-manifold, we define $M_{p,k}$ and $N_{p,k}$ as the result of some number of $(-1)$-torus surgeries and a $p$-torus surgery. This piece of information appears in the fifth column as a sum, where the first addendum refers to the number of $(-1)$-surgeries in the construction of $M_{p,k}$ and $N_{p,k}$ respectively. The bibliographic references containing the precise definition and the properties of the raw material and of the 2-tori involved in the surgical operations are listed in the last column. By construction, these families satisfy the hypothesis of Proposition \ref{Proposition S}, which proves the third item. It is straight-forward to verify that the mechanism of Section \ref{construction} applies to these examples and their equivariant intersection form is extended from the integers. The details are left to the interested reader.


\renewcommand{\arraystretch}{1.3}

\begin{table}\label{table2}
\begin{center}
\begin{tabular}{ |c|c|c|c|c|c| }
 \hline

 $k$ & $\chi$ & $\sigma$ & Raw material & \makecell{Number of\\[-3pt] surgeries} & Source \\
 \hline
 \hline

 $2$ & $5$ & $- 1$ & $(T^4\cs \CPbar^2)\cs _{\Sigma_2} (T^2\times \Sigma_2)$ & $(4 \text{ or } 3) + 1$ & \cite{[AkhmedovPark2], [Torres1]} \\ 
 
 \hline
 
 $3$ & $6$ & $-2$ & $(T^4\cs_2 \CPbar^2)\cs _{\Sigma_2} (T^2\times \Sigma_2)$ & $(4 \text{ or } 3) + 1$ & \cite{[AkhmedovPark1], [BaldridgeKirk], [Torres1]} \\ 
 
 \hline
 
$4$ & $7$ & $-3$ & $(T^4\cs \CPbar^2)\cs _{\Sigma_2} (T^4\cs_2 \CPbar^2)$ & $(2 \text{ or } 1) + 1$ & \cite{[AkhmedovPark2], [Torres1]} \\ 

 \hline
 
$5$ & $8$ & $-4$ & $(T^2\times S^2\cs_4 \CPbar^2)\cs _{\Sigma_2} (T^2\times \Sigma_2)$ & $(2 \text{ or } 1) + 1$ & \cite{[BaldridgeKirk], [Torres1]} \\ 

 \hline
 \end{tabular}
\end{center}
\vspace{3mm}
\caption{Raw materials for the collections with $k\in \{2, 3, 4, 5\}$}
\end{table}

\hfill $\square$





\end{document}